\documentclass[12pt]{article}

\usepackage{amscd}
\usepackage{amssymb,amsmath}

\newtheorem{theorem}{Theorem}
\newtheorem{prop}[theorem]{Proposition}
\newtheorem{lemma}[theorem]{Lemma}

\author{Vladimir Baranovsky}

\title{Uhlenbeck compactification as a functor.}
\date{Feb 23, 2014, revised Feb 10, 2015}

\begin{document}

\maketitle

\begin{center}
\thanks{Department of Mathematics, 
University of California - Irvine, \\Rowland Hall 340, 
Irvine, CA 92697 USA\\ email: vbaranov@math.uci.edu}
\end{center}

\section{Introduction.}
We work with noetherian schemes over an algebraically closed  field $k$ of arbitrary 
characteristic.  The goal of this paper is
to give a functorial definition of a compactification for moduli of
vector bundles on a smooth projective surface. 
More precisely, let $\pi: X \to S$ be 
a smooth projective morphism of relative dimension 2 
(although the theory can probably be generalized to the Cohen-Macaulay 
case). 
We define the \textit{Uhlenbeck moduli functor}, or the
functor of \textit{quasibundles} $QBun(r, d)$ that 
compactifies the stack $Bun(r, d)$ of vector bundles of fixed rank $r$
and $c_2 = d$ along the fibers of $\pi$.

It is well known, cf. \cite{DK}, 
that for $k = \mathbb{C}$, $S = Spec(\mathbb{C})$
one can use the Uhlenbeck 
compactness theorem to show that a family of vector bundles 
$E_t$, $u \in \mathbb{C} \setminus \{0\}$ with $c_2 = d$
degenerates,  in some sense, to a vector bundle $E_0$ with a different 
second Chern class, $c_2 (E_0) = d-k$ \textit{``plus"} an effective 
$0$-cycle $\xi$ of degree $k$ in $X$ , e.g. cf. 
\cite{DK}.  
In the original framework of the 
Uhlenbeck's Compactness Theorem, the cycle $\xi$ describes a 
non-negative combination of delta functions which split off the 
Chern form in the limit, for a family of connections on $E_t$. Understanding
the moduli space $M(r, d)$ in the sense of semistable bundles, one
gets a compact Uhlenbeck moduli space $M^U(r, d)$ which 
has a set theoretic decomposition
$$
M^U(r, d) = \coprod_{k \geq 0} M(r, d-k) \times Sym^k X
$$ 
and contains $M(r, d)$ as an open subset, although not necessarily
dense. Here $Sym^k$ is the $k$-th symmetric power of $X$,  i.e. the
quotient of the $k$-fold symmetric product by the action of the symmetric
group. In fact, by Bogomolov-Miyaoka inequality only finitely many pieces of this union are 
nonempty. 

In this paper, however, we would like work over an arbitrary 
algebraically closed field $k$ and consider the moduli
 functor  $QBun(r, d)$, containing rank $r$ 
vector bundles $Bun(r, d)$ as an open subfunctor, 
 with a certain compactness property 
(but non-separated - as is the case with $Bun(r, d)$
itself). 

A conceptual difficulty here is in treating 
formal sums $c_2(E) + \xi$ as a single object. 
The approach 
adopted in this paper is that both $c_2(E)$ and
 $\xi$ define 
multiplicative functors $Pic(X) \to Pic(S)$
between the groupoid categories of line bundles and their
isomorphisms. 

For $c_2(E)$ the corresponding functor is a 
categorical version of the map on cohomology
 $H^2(X, \mathbb{Z}) \to H^2(S, \mathbb{Z})$
or the Chow groups $A^1(X) \to A^1(S)$, obtained by taking
the cup product of $c_1$ for a line bundle $L$, with $c_2(E)$ and then integrating the class
along the fibers of  $\pi: X \to S$.  A fancier definition of 
the functor $\mathfrak{c}_2^E: Pic(X) \to Pic(S)$,  which the author
has learned from A. Beilinson, is that an $\mathcal{O}^*$-torsor $L^*$
on $X$ can be multiplied by the gerbe over the 
Zariski sheaf $\mathcal{K}_2$
of algebraic $K$-groups, corresponding to the class $c_2(E) \in H^2(X,
\mathcal{K}_2)$; and the resulting 2-gerbe over $\mathcal{K}_3$ may be
integrated along the fibers of $\pi$ to give a torsor over $\mathcal{O}^*$ on
$S$. Although in the last chapter we do make some comments on the gerbe
approach, the definition used in this paper for vector bundles
is that of $\mathfrak{c}_2^E(L)$ as the determinant of 
$$
R^\bullet \pi_*\  (L^\vee \to \mathcal{O}_X) \otimes_{\mathcal{O}_X} (\mathcal{O}_X^{r-1} \to E \to \det(E)).
$$
We can assume at first that both factors have zero arrows, but in fact
we can make them nontrivial by  choosing a section of $L$ and 
$(r-1)$ sections of $E$ - which is an important feature of this construction. The second
factor should be thought of as a categorical representative of $c_2(E)$ 
since the leading term of its Chern character is $c_2(E)$. 

 On the other hand, in characteristic zero
a family of zero cycles $\xi$ may be
viewed as a map $\xi: S \to Sym^d(X/S)$ to the relative symmetric power
of $X$. For fields of arbitrary characteristic
one needs to work with similar spaces $\Gamma^d(X/S)$ obtained via
the divided powers algebras, cf. \cite{Ro}, \cite{Ry}. 
For a line bundle $L$ on $X$ one has a natural line bundle $\Gamma^d(L)$
on $\Gamma^d(X/S)$ (informally, its fiber is the tensor product of the
fibers of the original $L$, over the $d$ points of the zero cycle).
Thus, for any section $\xi: S \to \Gamma^d(X/S)$
we can set $N_\xi(L) := \xi^* (\Gamma^d(L))$ which gives a functor 
$Pic(X) \to Pic(S)$.  

Both functors are\;\textit{multiplicative},
i.e. they agree with tensor products of line bundles in a natural 
sense. Tensor product of such functors will correspond 
to the formal sum $c_2(E_0) + \xi$. 

\bigskip
\noindent
The paper is organized as follows. First we define in Section 2 multiplicative
functors between Picard categories (in fact we deal with a functor 
$N: PIC(X/S) \to PIC(S)$ of larger categories that take
base changes into account); and give some examples.

In Section 3 we consider a quadruple $(Z, E, N, D)$ consisting of a
 vector bundle $E$ on  the open complement to  a closed  subset 
$Z \subset X$ which is finite over $S$, an extension of $\det(E)$
to a line bundle $D$ on $X$ and a multiplicative functor $N$ as before. 
In such a setting we define the notion of an 
$E$-localization of $N$ at $Z$  which says, informally, 
 that ``on the complement to $Z$ the functor $N$ is identified
with $\mathfrak{c}_2^E$". This is an additional structure, as the same 
$N$ can be localized at different closed subsets $Z$. Finally, in Section 4 we define
a notion of an ``effective $E$-localization", 
which boils down to saying that some apriori rational section of 
a line bundle is actually regular.  We need the concept to end up with 
a functor that eventually has a decomposition with pieces
like $M(r, d-k) \times Sym^k X$ before: without localization of functors
we would only be able to work with rational equivalence classes
of cycles, and without effectiveness - only with all cycles of fixed degree $k$
rather than effective ones.

In Section 5 we define the Uhlenbeck functor
as the functor of quadruples $(Z, E, N, D)$ equipped 
with an effective $E$-localization of $N$ at $Z$. Such 
objects will be called \textit{quasibundles} in this paper. 
Although the use of this term does not agree with \textit{some} earlier papers
in the field (see below) where ``quasibundles" stand for slightly different 
objects and the corresponding moduli spaces (functors) only 
map to the one defined here, we still choose this term due to 
the relation to quasimaps.  To explain it briefly, suppose that
we have a morphism of sheaves 
$\mathcal{O}^r \to E$
which is an isomorphism generically on each fiber over  $S$, and in fact on the 
complement to a relative curve $Y$ given by the determinant of 
$\mathcal{O}^r \to E$.  Then by Cramer's Rule
$E^\vee$ is sandwiched between two sheaves
$$
 \mathcal{O}_X^{r}(-Y) \subset E^\vee \subset \mathcal{O}_X^r
$$
so it can be described through a quotient of the sheaf $\mathcal{O}_Y^r$.
Further, let $(C_t, t \in \mathbb{P}^1)$ be a pencil of curves and
assume that none of the $C_t$ shares an irreducible component with $Y$. Then
for any $t$ we get a point in the Quot scheme $Quot_{C_t}
(\mathcal{O}_{Y \cap C_t}^r, l)$ parameterizing all 
quotients of $\mathcal{O}_{Y \cap C_t}^r$ with  fixed length $l$. 
Note that the scheme-theoretic intersection $Y \cap C_t$ is finite 
so the Quot scheme naturally embeds as a closed subscheme in a certain
product of Grassmannians, and its determinant bundle is the product of the
standard determinant bundles restricted from the Grassmannians. 
Varying $t$ we get a map (i.e. a section) from $\mathbb{P}^1$
to the corresponding relative Quot scheme, from which the original 
bundle $E$ can be recovered. 

For quasibundles this will be replaced by a quasimap in the sense of \cite{FGK}.
This single quasimap does not determine the quasibundle anymore: 
for a pair $(E, \xi)$ consisting of a bundle and a zero cycle, we can only recover
the direct image of $\xi$ to $\mathbb{P}^1$ (this makes sense if the base locus of the pencil 
avoids the support of $\xi$).  Still, working with a family of pencils, 
we can recover a quasibundle from the corresponding family
of quasimaps (this is explained in Section 4.5 to clarify the meaning of effective 
localization). 

\bigskip
\noindent
In Section 5.2 we show that
for  $S = Spec(k)$  the set of $k$-points of the stack of quasibundles has
the same decomposition as in the case of differential geometric
moduli space. We also prove that $QBun(r, n)$ is 
complete in a certain sense (i.e. establish the existence part in the 
valuative criterion of properness, but not the uniqueness which 
is not expected at all in this setting).  At the moment we cannot prove that the Uhlenbeck 
stack is algebraic but Section 5.3 gives a conjectured local flat covering by a scheme
of quasimaps into a relative punctual Quot scheme.  We also show that the Gieseker functor
of flat torsion free families of sheaves maps to the Uhlenbeck functor. 

\bigskip
\noindent
Finally, in Section 6 we explain how to define a
similar functor for torsors over split semisimple simply connected groups 
(the approach may be generalized to the reductive case). One important modification is
that effective zero cycles should be replaced by cycles with coefficients
in a semigroup. Then symmetric/divided powers of $X$ will be replaced by their
iterated products, as it should be if we expect that 
$$
QBun(G_1 \times G_2) = QBun(G_1) \times
QBun(G_2).
$$
In fact, we indicate how to define quasibundles in arbitrary dimension
but our treatment is necessarily sketchy
since this involves effective codimension 2 cycles, for which a functorial 
definition in arbitrary characteristic has not been written down yet. 

\bigskip
\noindent
We would like to place our work in regards to other papers on
compactifications for the moduli of bundles. The original topological
construction of Donaldson in \cite{DK} was reinterpreted in 
algebro-geometric terms independently by J. Li and J. Morgan, 
cf. \cite{Li} and \cite{Mo}.  In this approach the vector bundles 
are compactified by torsion free sheaves, and then two torsion
free sheaves $F$ and $G$ give the same point of the Uhlenbeck space if 
there is an isomorphism of vector bundles 
$F^{**} \simeq G^{**}$ and the two Artin sheaves $F^{**}/F$
and $G^{**}$ have the same support cycles. As was observed by a reviewer, 
this approach can also be implemented in finite characteristic. 
We also observe that the fact that  double duals are locally free will no longer 
hold when we consider families of bundles over a nontrivial base, and 
a instead of locally free or reflexive sheaves it is more convenient to work 
with those that satisfy a relative version of Serre's $S_2$ condition. Still, 
over a nilpotent base, one could expect that some families of quasi-bundles
in the sense of the present paper are not induced by families of torsion free 
sheaves (in a similar setting of the Hilbert-to-Chow morphism, there are 
families of effective cycles not induced by any flat family of Artin sheaves).

A framed version of this construction was implemented recently by 
U. Bruzzo, D. Markushevich and A. Tikhomirov in 
\cite{BMT}.

The general construction of V. Balaji
in \cite{Ba1} gives a compactified moduli scheme of principal $G$-bundles
of a surface over $\mathbb{C}$, 
which is obtained as a schematic image of another moduli scheme 
under a certain morphism. This scheme apriori depends on a choice
of a representation $\rho: G \to GL(V)$  but different choices gives
homeomorphic spaces.

In all of these cases the moduli scheme is defined as an image of another scheme
and it is not clear from the definition what kind
of objects it parameterizes. In a sense, this is the issue addressed in our
paper (with appropriate definition of semistability the moduli space of semistable
quasibundles - in our sense - is expected to be the scheme constructed in
\cite{Ba1}). In addition, the approach of the present paper works in 
arbitrary characteristic, and can be extended to arbitrary dimension. 
However, some aspects of this definition remain conjectural and
others need to be improved: 

\begin{itemize}

\item It would be very desirable to extend 
the constuction in Proposition 11 to the etale case and to establish its
degree property (see remark after the statement of this proposition). 

\item The ``categorical version" of $c_2$ for a principal bundle $P$
 used in this paper sends a line
bundle $L$ on $X$ to a line bundle on the base $S$, and we expect
that for a relative section $\sigma$ of $L$ this is isomorphic to the 
pullback of the canonical bundle on the moduli stack of $G$-bundles
over the family
of curves $Z(\sigma) \to S$.  In fact, it would be
much better to use this as a definition of our ``categorified" $c_2$
but for this one needs to extend the constuction of the canonical
bundle to the case of curves with nodal singularities (see remarks 
after the proof of Proposition 11).

\item Our definition of an effective localization is not intuitive and to 
improve it one would need a higher dimensional version of the 
Drinfeld-Simpson Lemma. 

\item Although we do obtain a functor $QBun(G, d)$
which is well behaved with 
respect to homomorphisms and products, and admits an 
expected decomposition involving (colored) zero cycles, cf. 
Proposition 12, it remains to be proved that an appropriate
representation $G \to \prod_i SL(r_i)$ induces a closed 
embedding, and that the usual $G$-bundles are open in $QBun(G, d)$.

\item Extension to the arbirary dimension depends on a similar functorial
construction for the Chow scheme of relative codimension 2 cycles
(work in progress by the author).  
\end{itemize}

In a different direction, a construction of the Uhlenbeck stack 
briefly outlined by V. Drinfeld in
\cite{Dr} has been implemented by A. Braverman, M. Finkelberg, 
D. Gaitsgory and A. Kuznetsov in \cite{FGK} and \cite{BFG} to yield 
an appropriate stack of framed quasibundles for the case of projective plane. 
In a sense, the present paper follows a similar approach mixed with
Deligne's theory of intersection bundles,
cf. \cite{De}, \cite{Du}, \cite{El}, \cite{MG}. It is rather likely that our 
definition reduces in the case of projective plane to one of the stacks defined in 
\cite{FGK} and \cite{BFG}.  Ideally, this approach should be phrased in 
terms of relative K-theory for flat surjective morphisms, but in this area there
has not been much progress since \cite{Gr1} and \cite{Gr2} (more specifically, 
one needs the analogue of global residue homomorphism in K-theory, and 
something like a relative Gersten resolution). See, however, \cite{Dr}, \cite{OZ} and a recent 
paper \cite{BGW} 
for the techniques on which such relative K-theory could be based. 
Perhaps the same range of methods can be used to define the 
categorical version of the second Chern class even in the non-projective
case, but then the language of gerbes seems to be unavoidable
(in our setting, still fairly abstract, we could get away with bundles due
to a version of the residue theorem). 

Other compactifications (or rather their stack versions) admit a natural 
``generalized blow-down" morphism to our functor of quasibundles.
These includes the Gieseker compactification by torsion free sheaves, 
as long as compactifications by G\'omez-Sols in \cite{GS} and A. Schmitt
in \cite{Sch}, the algebraic version of the bubble-tree compactification by 
D. Markushevich, A. Tikhomirov and G. Trautmann in \cite{MTT} 
and a somewhat similar compactification by N. Timofeeva in \cite{Ti}. 
All of these compactifications remember more than just the multiplicities
of the cycle $\xi$ and consequently none of such moduli functors satisfies
the product formula 
$\mathcal{M}(G_1 \times G_2) \simeq \mathcal{M}(G_1) \times 
\mathcal{M}(G_2)$.

Finally, a word on terminology. Our use of the word ``quasibundles"
agrees with that of \cite{FGK} and \cite{BFG} (since ultimately they are
linked to quasimaps) and differs from that of \cite{Ba1}, \cite{Sch} 
and \cite{GS}. 

All unlabeled tensor products are tensor products of quasicoherent sheaves
on a scheme over its structure sheaf. 

\bigskip
\noindent
\textbf{Acknowledgements.} The author has
greatly benefitted from conversations with A. Beilinson, 
V. Drinfeld, V. Ginzburg, A. Kuznetsov and M. Finkelberg, and 
expresses his deep gratitude to them. The author also thanks the reviewers 
for the suggested improvements.

\section{Multiplicative functors}

For a morphism $\pi: X \to S$ we denote by $PIC(X/S)$  the category
 formed by  pairs $(T, L)$ where $T$ is a scheme over $S$ and $L$ is a line bundle
over $X_T = X\times_ST$. A morphism $(T_1, L_1)\to 
(T_2, L_2)$ is an $S$-morphism $T_1 \to T_2$ plus an isomorphism
of $L_1$ with the pullback of $L_2$ to $X_{T_1}$. 
When $X = S$ and $\pi$ is the identity map, we denote $PIC(X/S)$ by
$PIC(S)$.

We choose and fix a locally constant function $d: S \to \mathbb{Z}$,
such as the second Chern class for a family of bundles with base $S$
(when $\pi$ is a family of smooth projective surfaces); or the 
degree of a family of zero cycles with base $S$.

\medskip
\noindent
\textbf{Definition.} By a degree $d$\;\textit{multiplicative functor} we will 
understand the data consisting of 

(a) A functor $N: PIC(X/S) \to PIC(S)$,

(b) A family of isomorphisms $b_{L_1, L_2}: N(L_1 \otimes L_2) 
\to N(L_1) \otimes N(L_2)$ which agree with the associativity and commutativity
for tensor products of line bundles; in the sense of
Section 1.2.1 of \cite{Du}.

(c) A family of isomorphisms $c_{(T, L)}: N(\pi^*_T L) \simeq L^{\otimes d}$, 
where $L$ is a line bundle on $T \to S$ and $\pi_T: X_T \to T$ is
the base change. This should agree with isomorphism of (b) 
in a natural way, and we also require that
 an isomorphism $\sigma: L_1 \to L_2$ corresponds to 
$c_{\sigma} =\sigma^{\otimes d}: L_1^{\otimes d} \to L_2^{\otimes d}$. 

\bigskip
\noindent
We now give three basic examples of multiplicative functors

\subsection{Zero Cycles and Families of Artin Sheaves}

In characteristic zero, let $Sym^d(X/S) \to S$ be the relative symmetric power
obtained by taking the quotient of the $d$-fold relative cartesian product of $X$ over
$S$, by the action of the symmetric group $\Sigma_d$. If $\pi$ is projective,
so is $Sym^d(X/S) \to S$. In characteristic $p> 0$ one needs to be more careful
and work with divided powers, cf. \cite{Ro}, \cite{Ry}, 
which leads to a similar space
$\Gamma^d(X/S) \to S$. For any line bundle $L$ on $X$ one has an induced
line bundle $\Gamma^d(L)$ on $\Gamma^d(X/S)$;  in characteristic 
zero this is simply the descent of $L \boxtimes \ldots \boxtimes L$ with 
respect to the symmetric group action. 

A family of zero cycles of degree $d$ with base $S$ is represented by 
a section $\xi: S \to \Gamma^d(X/S)$ and therefore we have a 
functor $N: PIC(X/S) \to PIC(S)$ which sends $(T \to S, L)$ to a line
bundle $\xi_T^* (\Gamma^d L)$ on $T$. Here $\xi_T: T \to 
\Gamma^d(X/S)_T \simeq \Gamma^d(X_T/T)$ is the family 
of zero cycles obtained by base change. 
It is not too difficult to construct the additional data required in the 
definition of a multiplicative functor, see \cite{Ba}.

\bigskip
\noindent
A family of zero cycles can be induced by a sheaf $\mathcal{A}$  
on $X$ which is flat with finite support over $S$ (i.e. a family 
of Artin sheaves with base $S$). In this case $\xi$ simply gives 
points in the support of $\mathcal{A}$ with 
multiplicities. Then the multiplicative functor can also be defined by
$$
N_{\mathcal{A}} (L) = \det \pi_* \big[(L^\vee \to \mathcal{O}_X) \otimes \mathcal{A}\big]
$$
See  Chapter 6 of \cite{De}, Section 1.2.1 of \cite{Du} for more 
details. By Section 5 in \cite{Ba} this is the functorial version of the 
Hilbert-to-Chow map in the case when $\mathcal{A}$ is a family
of zero dimensional subschemes.

\subsection{Intersection Bundles}

Another example of a multiplicative functor $PIC(X/S) \to PIC(S)$
is the two-dimensional version of Deligne's intersection bundles
introduced in \cite{De}. For
three bundles $L_0, L_1, L_2$ on $X$ we set
$$
\langle L_0, L_2, L_2 \rangle:=
\det R\pi_* \big[(L_0^\vee \to \mathcal{O}_X) \otimes 
(L_1^\vee \to \mathcal{O}_X) \otimes 
(L_2^\vee \to \mathcal{O}_X)\big]
$$
The properties of these construction were studied extensively in 
\cite{El}, \cite{MG}, \cite{Du}. 
Fixing $L_1$ and $L_2$ one obtains a multiplicative functor
$$
IB_{L_1, L_2}: PIC(X/S) \to PIC(S); \qquad L_0 \mapsto 
\langle L_0, L_1, L_2 \rangle.
$$

\subsection{Second Chern Class}

In this subsection we give two different constructions 
for the second Chern class functor $\mathfrak{c}_2^E:
PIC(X/S) \to PIC(S)$ of a bundle $E$ on $X$, and briefly 
explain why they are equivalent. The second construction is 
a special case of a more general definition given earlier by 
R. Elkik; and the first one (which seems to be new) is a bit more 
adjusted to the concept of a localized multiplicative
functor, that we introduce below. 

\bigskip
\noindent
\textbf{Construction 1.} The starting point here is a geometric 
interpretation of the second Chern class on a smooth projective
$X$.
If $E$ is generated by global sections then the rank of a
generic morphism $\xi: \mathcal{O}_{X}^{r-1} \to E$
will be $\leq (r-2)$ on a codimension 2 closed subset  $Z_\xi$ 
which corresponds to the class $c_2(E)$ in an appropriate 
version of cohomology group (e.g. Chow groups or integral 
cohomology if $k = \mathbb{C}$) - of course, if one takes into 
account multiplicities. 
We can restate this using a particular case of the
 \textit{Buchsbaum-Rim
complex}, cf. \cite{BR}. 
Let $\xi: F \to E$ be a morphism of vector bundles of
ranks $(r-1)$ and $r$, respectively. Consider the  \textit{dual}
Buchsbaum-Rim
complex $\mathcal{BR}_\xi$ 
\begin{equation}
\label{KofE}
0 \to F \to E \to \det E \otimes \det F^\vee \to  0
\end{equation}
where the first nonzero arrow is  $\xi$, and the second is 
obtained from the map $E \otimes \Lambda^{r-1} F \to \Lambda^r E$
sending $e \otimes \omega$ to $e \wedge \Lambda^{r-1}(\xi)(\omega)$.
We will assume that the term $\det E \otimes \det F^\vee$ 
is placed in homological degree zero. A simple argument 
shows that
this complex is exact away from $Z_\xi$ and has Chern character 
\begin{equation}
ch(\mathcal{BR}_\xi) = c_2(E) - c_2(F) + c_1(F)^2 - c_1(F) c_1(E) +
 \textrm{higher\;order\;terms}
\end{equation}
When $F = \mathcal{O}_X^{r-1}$ this reduces to 
$\big(c_2(E) +$ higher order terms$\big)$.
\begin{prop}
Suppose that $\pi: X \to S$ is a flat family of Cohen Macaulay
projective surfaces and the degeneracy locus $Z_\xi$ of
$\xi: \mathcal{O}^{r-1}_X \to E$ is 
finite over $S$. Then $\mathcal{BR}_\xi$ has non-trivial 
cohomology in degree zero only, and  the
sheaf $\mathcal{H}^0(\mathcal{BR}_\xi)$ is flat over $S$. 
\end{prop}
\textit{Proof.} By 
definition $\mathcal{BR}_\xi$ is dual to the Buchsbaum-Rim
complex, cf. \cite{BR}, 
of the adjoint morphism $\xi^*: E^\vee\to  \mathcal{O}^{r-1}_X$. 
Since $\pi$ is Cohen-Macalay and $Z_\xi$ is finite, 
for any point $x \in Z_\xi$ and any sheaf $\mathcal{E}$ on $S$, 
there exists a length two $\pi^* \mathcal{E}$-regular sequence
in $\mathfrak{m}_x$. By Theorem 1 of \cite{BR}, 
for any such $\mathcal{E}$ the complex 
$\mathcal{BR}_\xi \otimes \pi^* \mathcal{E}$ can 
only have cohomology in degree zero. 
For $\mathcal{E} = \mathcal{O}_S$ we
obtain the statement about the cohomology of $\mathcal{BR}_\xi$.
Viewing $\mathcal{BR}_\xi$ as a resolution of 
$\mathcal{H}^0(\mathcal{BR}_\xi)$ 
we conclude also that 
$\mathcal{T}or_{i}^{\mathcal{O}_S}(\mathcal{H}^0(\mathcal{BR}_\xi), 
\mathcal{E}) = 0$ for $i > 0$, i.e. $\mathcal{H}^0(\mathcal{BR}_\xi)$ is 
flat over $S$. $\square$

\medskip
\noindent
We can now define a functor
$$
\mathfrak{c}_2^E: PIC(X/S) \to PIC(S)
$$
as the functor $N_{\mathcal{A}}$ defined for the family 
of Artin sheaves $\mathcal{A} = \mathcal{H}^0(\mathcal{BR}_\xi)$
as in Section 2.1.
To show that it depends only on $E$ but not $\xi$ and also to 
define it for bundles which are not generated by global sections,
observe that by 
 the Euler isomorphism, cf. e.g. \cite{GKZ},
  $N_{\mathcal{A}}$ is isomorphic to 
the
functor 
\begin{equation}
\label{defineC2}
L \mapsto \det\; R \pi_{S*} \big[(L^\vee \to \mathcal{O}_{X})\otimes
(\mathcal{O}^{r-1}_{X} \to E \to \det E)\big]
\end{equation}
where the arrows in the two tensor factors are zero. 

\bigskip
\noindent
\textbf{Remark.}
More generally, it is natural to consider a situation when we have
a  rank $(r-1)$ vector bundle 
$G$ on $S$ and a coherent sheaf morphism $\xi: W=\pi^* G \to E$
which is an embedding of vector bundles 
away from a closed subset $Z_\xi$, finite over $S$.
This leads to a complex 
$$
0 \to \pi^* G \to E\to \det E \otimes \det \pi^*G^\vee \to 0
$$
An easy application of the projection formula shows that 
$$
\det R\pi_* \big[(L^\vee \to \mathcal{O}_X) \otimes
\mathcal{BR}_{\xi} \big]\simeq \mathfrak{c}_2^E(L) \otimes
\det G^{- L \cdot \det (E)};
$$ 
where for two line bundles $L, D$ we denote 
$$
L \cdot D = \chi (L \otimes D) - \chi(L) - \chi(D) + 
\chi(\mathcal{O}_X)
$$
and $\chi$ is the Euler characteristic of a vector bundle restricted to 
fibers of $\pi$. By flatness of $\pi$ this is a locally constant 
function $S \to \mathbb{Z}$.

\bigskip
\noindent
\textbf{Construction 2.} A general construction of
R. Elkik in \cite{El} specializes to our case 
as follows. First consider the projective bundle
$\psi: \mathbb{P}:=\mathbb{P}(E^\vee) \to X$ of lines in $E$, 
and the corresponding Segre class functor $\mathfrak{s}_2^E:
PIC(X/S) \to PIC(S)$ given by 
$$
\mathfrak{s}_2^E(L) = \det R\psi_* \big[(\psi^* L^\vee \to \mathcal{O}_{\mathbb{P}}) 
\otimes (\mathcal{O}_{\mathbb{P}}(-1) \to 
\mathcal{O}_{\mathbb{P}})^{\otimes (r+1)} \big].
$$
Next, imitating the
relation between Chern and Segre classes $c_2(E) + s_2(E) = c_1^2(E)$ we
set
$$
\mathfrak{c}_2^E(L):= IB_{\det(E), \det(E)}(L) 
\otimes \mathfrak{s}_2^E(L)^\vee  
$$
Let us explain briefly why this
functor is isomorphic to that of Construction 1. First, one uses the 
fact that $\mathfrak{s}_2^E$ is multiplicative, 
cf. \cite{El},  to show that the second definition 
is equivalent to
$$
L\mapsto \det R(\pi \psi)_* \big[(\psi^* L^\vee \to
\mathcal{O}_{\mathbb{P}}) \otimes (\mathcal{O}_{\mathbb{P}} (-1)\to 
\mathcal{O}_{\mathbb{P}})^{\otimes r} \otimes 
(\mathcal{O}_{\mathbb{P}} \to \mathcal{O}_{\mathbb{P}}(1))\big].
$$
The projection formula for the morphism $\psi$ shows that
$\mathfrak{s}_2^E$ is isomorphic to the functor
$$
L \mapsto \det 
R \pi_* \big[(L^\vee \to \mathcal{O}_X) \otimes
(\det(E)^\vee \to \mathcal{O}^{r+1} \to E) \big]
$$
(as before, all arrows in the factors on the right hand side are zero). Now the 
isomorphism of the two definitions of $\mathfrak{c}_2^E$ 
follows from a similar isomorphism
$$
IB_{\det(E), \det(E)} (L) 
\simeq \det R \pi_* \big[(L^\vee \to \mathcal{O}_X) 
\otimes (\det(E)^\vee \to \mathcal{O}_X) \otimes
(\mathcal{O}_X \to \det(E))\big]
$$
$$
\simeq \det R\pi_* (L^\vee \to \mathcal{O}_X)\otimes 
(\det(E)^\vee \to \mathcal{O} \oplus \mathcal{O} \to \det(E))
$$
Thus, the streamlined version 
of Construction 2, is given by:
$$
\mathfrak{c}_2^E(L):= IB_{\det(E), \det(E)}(L) 
\otimes \mathfrak{s}_2^E(L)^\vee  
$$
where 
$$
\mathfrak{s}_2^E(L) \simeq \det R \pi_*
(L^\vee \to \mathcal{O}_X) \otimes (\det(E)^\vee \to \mathcal{O}^{r+1} \to E)
$$
Comparing the definitions we see that the above formula
gives the same functor as Construction 1.

\section{Localization of multiplicative functors}

Let $Z \subset X$ be a closed subset, finite over $S$. Consider a bundle $E$ 
on the open subscheme $U = X \setminus Z$.
To be able to say that a multiplicative functor $N: PIC(X/S) \to PIC(S)$
``restricts to $\mathfrak{c}_2^E$ over $U$" we need the concept of an
\textit{$E$-localization of $N$ at $Z$}. By definition, this is the additional data 
described as follows. 

Consider a morphism $\phi: L_1 \to L_2$ of line bundles on $X$ viewed as
coherent sheaves, and assume that $\phi$ is an isomorphism along $Z$. 
Then $Coker(\phi) \otimes (\mathcal{O}^{r-1} \to E \to \det(E))$ is 
a well-defined perfect complex on $X$: although the second factor is only defined
on $U$, the support of the first factor is contained in $U$.  To simplify notation
we will assume $L_2 \simeq \mathcal{O}_X$ but the general case can 
be derived easily. 

\medskip
\noindent
\textbf{Definition.} An 
$E$-localization of $N$ at $Z$ consists of a family of isomorphisms 
$$
a_{\phi}:  \det \pi_* \big[\mathcal{O}_Y \otimes
 (\mathcal{O}^{r-1} \to E \to \det(E))\big] \simeq N(L),
$$ 
for any $\phi: L^\vee \to \mathcal{O}$ which vanishes on 
a relative divisor $Y$ disjoint from $Z$. We require that
such isomorphisms should agree with tensor product on sections 
and values of $N$. 

\bigskip
\noindent
We have the following basic examples

\medskip
\noindent
\textbf{Supports of zero cycles.} If $N: PIC(X/S) \to S$ is the functor 
describing the family of zero cycles, then $N$ carries a canonical 
$\mathcal{O}$-localization
at the closed subset $Z_\xi$ which is the support of the relative zero cycle (or the 
support of the family of Artin sheaves). 
We remark here that a rigorous definition of the support for a family
$\xi: S \to \Gamma^d(X/S)$ of zero cycles, involves
the addition map $add\!: X \times_S \Gamma^{d-1}(X/S) \to \Gamma^d(X/S)$, cf. \cite{Ry}, \cite{Ba}.
Then the support $Z_\xi$ is simply the projection of the closed subset
$add^{-1} (\xi(S))$ to $X$. 

\medskip
\noindent
\textbf{Second Chern class.} For a bundle $E$ on $X$  the definition of
$N = \mathfrak{c}_2^E$ implies that this multiplicative functor has a
canonical $E$-localization at $Z = \emptyset$.

\medskip
\noindent
\textbf{Generic morphism $\mathcal{O}^{r-1} \to E$.} Suppose we have 
a morphism $\phi: \mathcal{O}^{r-1} \to E$ which has rank $\leq (r-2)$ 
on a closed subset $Z \subset X$. Then the non-zero arrows of the dual
Buchsbaum-Rim complex
$\mathcal{O}^{r-1} \to E \to \det(E)$ make it quasiisomorphic to a
family of Artin sheaves $\mathcal{A}$ supported at $Z$, hence this induces
an $\mathcal{O}$-localization of $\mathfrak{c}_2^E$ at $Z$. 

\bigskip
\noindent
The last example illustrates that an $E$-localization is really some additional 
data which cannot be recovered from $N$ alone - different choices of 
$\phi$ will in general lead to different closed subsets $Z$. Hence the 
same functor $N$ may admit localizations 
at different closed subsets. This is somewhat similar to the fact that an effective divisor
on a curve defines a line bundle but different effective divisors may give isomorphic
line bundles.  

\begin{lemma} Any functor automorphism of $c$, which agrees with
the multiplicative structure and the
localization on $Z$, is trivial.
\end{lemma}
\textit{Proof.} By definition, a functor automorphism of
$c$ should send a pair $(L, T)$ to an automorphism
of the line bundle $c(L)$ on $T$. Agreement with 
the multiplicative structure of $c$
 means that for $T \to S$ we have a 
group homomorphism  $\eta: Pic(X_T) \to \Gamma(\mathcal{O}^*_T)$, from the 
group $Pic(X_T)$ of isomorphism classes of
line bundles on $X_T$, 
which behaves naturally with respect to pullback
on $S$-schemes.

To show that $\eta$ is trivial on the class of $L$
we may assume that $T$ is affine and $L$ is very ample.
Then we can cover $T$ with smaller open 
subsets $T_i$ and choose over each $T_i$ a
section $\phi_i: \mathcal{O}_{X_T} \to L$ which 
vanishes away from $Z_{T_i}$.  By agreement with 
localization, over each $T_i$ we must have 
$\eta(L) = \eta(\mathcal{O}_{X_T}) = 1$, as 
required. 
$\square$

\section{Effectiveness conditions for localizations}

Effectiveness is a condition imposed on localizations, and roughly speaking
it means that some apriori rational section of a line bundle is actually regular. We first 
explain how this condition appears in examples, where it is based on the 
following result, obtained by combining  Theorem 3, Proposition 9 and
Lemma 1 on p. 51 in \cite{MK} (the notation and 
formulation is changed to match our situation):

\begin{prop}\label{mk}
Let $\pi: X \to Y$ be a proper morphism of noetherian 
schemes of finite Tor-dimension, and assume that 
$Y$ is normal. If 
$\mathcal{F}^\bullet$ is a perfect complex on $X$ such that
\begin{enumerate}
\item for each depth 0 point $y \in Y$,  restriction of $\mathcal{F}^\bullet$
to the fiber over $y$ is exact, 
\item for each depth 1 point $y \in Y$,  restriction on $\mathcal{F}^\bullet$
to the fiber over $y$ has cohomology only in degree zero, and this unique
cohomology sheaf has $0$-dimensional support. 
\end{enumerate}
Then the natural trivialization of $\det R \pi_*(\mathcal{F})^\bullet$ over an 
open subscheme of $Y$ (which exists due to (1))
extends canonically to a regular section 
$Div(\mathcal{F}^\bullet)$
of the line bundle  $\det R \pi_* (\mathcal{F}^\bullet)$ and for any line bundle
$L$ on $X$ one has $Div(\mathcal{F}^\bullet) = Div(\mathcal{F}^\bullet \otimes L)$.
\end{prop}

\subsection{Zero cycles and generic sections of  vector bundles.}

Suppose we are in a situation of Section 2.1 when there is 
a flat family of Artin sheaves $\mathcal{A}$ on $X$ and assume for 
simplicity in this subsection that $S = Spec(k)$. Choose a  line bundle
$L$ on $X$ generated by a subspace $H$ of its global sections
and set $P = \mathbb{P}(H)$. Then the bundle 
$\mathcal{O}_P(1) \boxtimes L$ on $X_P = P\times_k X$ has a 
canonical section $s$. If the arrow in the 
complex $(\mathcal{O}_P(-1) \boxtimes \mathcal{A}(L^\vee)
\to \mathcal{O}_P \boxtimes \mathcal{A})$ on $X_P$
is given by $s^\vee \otimes
Id_{\mathcal{A}}$ then this complex becomes exact over the 
open subset $U_P \subset P$ of all curves in the linear system which 
avoid the finite support of $\mathcal{A}$.  The second condition 
of the previous Proposition is also satisfied. By the result quoted above
the trivialization of $\det R \pi^* \big[\mathcal{O}_P(-1) \boxtimes \mathcal{A}(L^\vee)
\to \mathcal{O}_P \boxtimes \mathcal{A}\big] \simeq \mathcal{O}_P(d)$
over the open subset $U$ extends to a regular section $Div$. In fact, each point
$x \in Supp(\mathcal{A})$ corresponds to a hyperplane $H_x \subset Y$
of curves in the linear system that contain $x$, and the section $Div$
is the product of  the linear equations for $H_x$ each repeated 
$m_x$ times, where $m_x$ is the multiplicity of $x$ in the support cycle
of $\mathcal{A}$. 

\bigskip
\noindent
Now suppose that a morphism $\psi: \mathcal{O}_X^{r-1} \to E$ which 
has rank $\leq (r-2)$ at a finite subset $Z$. Then the functor $\mathfrak{c}_2^E$ 
acquires an $E$-localization at $Z$: we can first replace the zero arrows in 
$\mathcal{O}_X^{r-1} \to E \to \det(E)$ by the arrows of the dual 
Buchsbaum-Rim complex and then replace the complex itself by its cohomology
sheaf $\mathcal{A}$ in degree zero. Again, for a line bundle $L$ the trivialization of 
$\mathcal{O}_P(d)$ at the open subset of sections vanishing away from $Z$, extends
to a regular section. 

We also observe here that if  $E$ is generated by global sections then a sufficiently 
general $\varphi$ will have rank $\leq (r-2)$ at most at a finite subset of $X$. 
Indeed, choose a surjection $W \otimes_k \mathcal{O}_X \to E$ and 
let $G = Gr(r-1, W)$.  If $\mathcal{Z} \subset G \times X$ is the closed subset
of pairs $(\psi, x)$ such that $\psi$ has rank $\leq (r-2)$ at $x \in X$. 
The for every $x \in X$ its preimage in $G$ is a standard Schubert variety of
codimension 2. Therefore, $\mathcal{Z}$ is irreducible of the same dimension as $G$.
Hence for the projection $\mathcal{Z} \to G$, a generic $\varphi \in G$  has
a finite fiber, as required.

\bigskip
\noindent
\textbf{Remark.} 
For a line bundle $L$, let $H$ be a vector subspace of $H^0(X, L)$ which generates
$L$ and $W$ as above. Using similar arguments one can 
construct a similar trivialization of the line bundle
$\mathcal{O}_P(d) \boxtimes \mathcal{O}_G (L \cdot \det(E))$ 
 on the open subset of all pairs  $(s, \psi) \in P \times G$ for which 
$\psi$ has maximal rank $(r-1)$ everywhere at the zero set of $s$; 
and also extend it to a \textit{regular} section of 
$\mathcal{O}_P(d) \boxtimes \mathcal{O}_G (L \cdot \det(E))$. This
generalizes the classical Chow form for cycles in a projective space. 

\subsection{Twists by a line bundle.}

A general bundle $E$ will be generated  by a finite dimensional space of its 
sections $W$ only after a twist by an appropriate very ample line bundle
(in a relative situation $W$ would be a pullback for some bundle on $S$). 
We would like to explain how the functor $\mathfrak{c}_2$ changes
when $E$ is replaced by $E(M)$ for some line bundle $M$. On the level of 
cohomology classes, if $x_i$ are the Chern roots of $E$ and $h$ is the 
first Chern class of $M$ then $x_i + h$ are the Chern roots of $E(M)$.
Therefore we have a formula
$$
c_2(E(M)) = c_2(E) + (r-1) h c_1(E) + \frac{r(r-1)}{2} h^2.
$$
Its categorical version is an isomorphism of functors
$$
\mathfrak{c}_2^{E(M)} \simeq \mathfrak{c}_2^E \otimes IB_{M, \det(E)}^{\otimes (r-1)}
\otimes IB_{M, M}^{\otimes \binom{r}{2}};
$$
and its proof is an easy exercise left to the interested reader (see the end of \cite{MK} for a general 
approach). Denoting by $c(\det(E), M)$ the last two factors, we
can write it in the form
$$
\mathfrak{c}_2^{E(M)} \simeq \mathfrak{c}_2^E \otimes c(\det(E), M) .
$$
Similarly, for any multiplicative functor $N: PIC(X/S) \to PIC(S)$
we will write 
$$
N^{(M)} = N \otimes c(\det(E), M).
$$
It is straightforward to check that the twisted cycle has degree
$$
d(M) = d + (r-1) M \cdot \det(E) + \binom{r}{2} M \cdot M.
$$
It follows from the definitions that an $E$-localization of $N$ at $Z$ 
induces an $E(M)$ localization of $N(M)$ at $Z$.

\subsection{Effective localizations.}

Now we return back to the situation when $E$ is defined only on the 
complement $U$ of a closed $Z \subset X$ which is finite along the fibers 
over  $S$. 
We fix a multiplicative functor $N: PIC(X/S) \to PIC(S)$ and its
$E$-localization at $Z$, which means that for a section   
$\phi: \mathcal{O}_X \to L$ 
vanishing on a relative curve $C$ disjoint from  $Z$, we 
are given an isomorphism
$$
a_{\phi}:  \det R\pi_* \big[\mathcal{O}_C \otimes
 (\mathcal{O}^{r-1} \to E \to \det(E))\big] \simeq N(L),
$$ 
where $\phi^\vee: L^\vee \to \mathcal{O}_X$ is the adjoint.

\bigskip
\noindent
After a base change $T \to S$  and a twist by an ample line bundle $M$
on $X_T = X \times_S T$ we can find a morphism 
$\psi: \mathcal{O}^{\oplus(r-1)} \to E_T(M)$ on the open complement 
$U_T$ to the closed subset $Z_T$, which has maximal rank $(r-1)$
away from a closed subset of $U_T$ which is also finite over $T$. 
Denoting by $Z_\psi$ the union of this subset with $Z_T$ we have that
$\psi$ has maximal rank $(r-1)$ on the open complement $U_\psi$ to 
$Z_\psi$. But then we can choose nonzero arrows in 
$(\mathcal{O}^{r-1} \to E_T(M) \to \det(E_T(M)))$ as in the dual Buchsbaum-Rim
complex of $\psi$ and obtain an exact complex on $U_\psi$. Thus, an $E$-localization
of $N$ and a choice of $\psi$ induce a trivialization of $N^{(M)}(L)$ for any
section $\phi$ of $L$ which is nonzero along $Z_\psi$, and this trivialization agrees
with tensor products on line bundles an their sections. 

In other words, we have extended the domain of  $N^{(M)}$ to the 
category $PIC_{Z_\psi}(X_T/T)$ which still has line bundles as objects but 
now morphisms are coherent sheaf morphisms which are nonzero at every point of 
$Z_\psi$.
By \cite{Ba2} such a functor must come from a relative zero cycle on $X_T$, but
not necessarily effective (as a difference of two 
effective cycles will also be trivially localized at the union of 
their suppors). This motivates the following definition:

\bigskip
\noindent
\textbf{Definition.} An $E$-localization $a$ of $N$ is \textit{effective} if 
for any base change $T \to S$ and a morphism $\psi: \mathcal{O}^{r-1} 
\to E_T(M)$ such that $M$ is relatively ample and $\psi$ has maximal rank 
on  the complement of $Z_\psi \subset X_T$
as above, the functor $N^{(M)}$ with its $\mathcal{O}$-localization at
$Z_\psi$ is isomorphic to the functor of an effective zero cycle 
$\xi: T \to \Gamma^{d(M)} (X_T/T)$ with support in $Z_\psi$.

\bigskip
\noindent
For the readers convenience we recall that this means that the extension 
$$
N^{(M)}: PIC_{Z_\psi}(X_T/T) \to PIC(T)
$$ 
actually 
extends further to a functor $PIC^+(X_T/T) \to PIC^+(T)$ 
where both $PIC^+$ have the same objects as before, but now
morphisms are arbitrary morphisms of coherent sheaves. This is actually a certain 
condition imposed on the extension to $PIC_{Z_\psi}$: if we have a family 
of sections $\phi_t$, $t \in \mathbb{A}^1$ for a line bundle $L$ and 
 for $t \neq 0 \in \mathbb{A}^1$ the zero set of $\phi_t$ avoids $Z_{\psi}$,
then the pullback of $N(L)$ to $\mathbb{A}^1 \times_k T$ is trivialized on 
$(\mathbb{A}^1 \setminus 0) \times_k T$ and we require that this trivialization 
extends to a regular section on $\mathbb{A}^1 \times_k T$.

\bigskip
\noindent
Below we interpret the effectiveness condition in terms of quasimaps to punctual 
Quot schemes of curves in $X$.

\subsection{Relation to effective zero cycles.}

\begin{prop}
Suppose that $E$ is a bundle defined everywhere on $X$ and $N: PIC(X/S)
\to PIC(S)$ is a multiplicative functor with an effective $E$-localization 
at a closed subset $Z \subset X$ finite over $S$. Then there exists
an effective $0$-cycle $\xi$ of degree 
$k = \deg N - \deg \mathfrak{c}_2^E$ with support on $Z$ and a 
multiplicative isomorphism of functors that agrees with localizations
$$
N \simeq \mathfrak{c}_2^E \otimes N_\xi
$$

\end{prop}
\textit{Proof.} By assumption $\mathfrak{c}_2^E$ is well defined 
so we can form the multiplicative functor $N_1 = N \otimes (\mathfrak{c}_2^E)^{-1}$.
Comparing $E$-localizations we see 
that for any section $\phi: \mathcal{O}_X \to L$ of a line bundle
which does not vanish at $Z$ (perhaps chosed after a base change
on $S$), we are given an isomorphism
$$
a_{\phi}: \mathcal{O} \simeq N_1(L);
$$
and such isomorphisms behave multiplicatively with respect to 
tensor products of line bundles and their sections. 

We would like to show that $a_{\phi}$ can be extended to 
a coherent sheaf morphism even in the case when $\phi$ does vanish
at the points of $Z$.  Since $N_1$ does not change after twist by a
relatively ample $M$, and effectiveness can be detected by values on
sufficiently ample $L$, we can localize on $S$ and 
assume that $L$ and $E$ are generated by  
finite dimensional spaces $H$, $W$ of their sections.  
Set $P = \mathbb{P}(H)$ as
before.

Over the open subset $\mathcal{U} \subset S \times_k P$ of sections that do not 
vanish at $Z$ the bundle $N_1(L \boxtimes \mathcal{O}_P(1))
= N_1(L) \otimes \mathcal{O}_P(k)$ is trivialized by assumption. 
We need to explain why the trivializing section $a_\phi$ is regular
along  $P \setminus \mathcal{U} = Z_P$ (fiber by fiber over $S$ our $Z_P$ is
a union of hyperplanes in $P$ formed by curves passing through 
points of $Z$). To that end, we can
assume $S = Spec(k)$. 
Choose a point $p \in Z_P$, which corresponds to a curve 
$C_p \subset X$.  In the Grassmannian $G$ of $(r-1)$-dimensional
subspaces in $W$, the subscheme of all $g \in G$ which induce
a linear map $U_g \to W \to E_x$ of non-maximal rank $\leq (r-2)$
is a codimension 2 Schubert cell, for any $x \in C_p$. Therefore we 
can choose a subspace of sections $U_g$ of dimension $(r-1)$ which 
has maximal rank everywhere on $C_p$.  Then, in the neighborhood
of $p$ we obtain a trivialization of $\mathfrak{c}_2^E(L \otimes 
\mathcal{O}_P(1))$ and its product with $a_\phi$ is a regular
section of $N(L \otimes \mathcal{O}_P(1))$, by assumption. 
Hence $a_\phi$ itself was regular in the neigborhood of $p$ to
begin with.  

Therefore, the effective localization of $N$ extends $N_1$ to 
 a multiplicative degree $k$ functor 
$PIC(X/S)^+ \to PIC(S)^+$ of larger categories which still have 
line bundles as objects but morphisms are taken in the sense of 
coherent sheaves. By the main result of \cite{Ba},  such a functor
is isomorphic to a functor $N_\xi$ associated to a unique 
degree $k$ effective cycle $\xi$. 
$\square$

\bigskip
\noindent
Perhaps this is an appropriate place to explain why just
considering bundles $E$ away from finite subsets is not
good enough, and also a functor $N$
with an effective localization is needed. If $X$ is a smooth surface over
a field then any vector bundle $E$ on a complement $U$ to 
a finite set $Z \subset X$ extends uniquely to $X$: its direct image
as a coherent sheaf satisfies Serre's $S_2$ condition
 and then by Auslander-Buchsbaum formula the 
direct image should be projective at the points of $Z$. However, 
the second Chern class of this extension does not remain constant
when we start varying $E$ and $Z$ in a family. Yet we 
have a semi-continuity statement. 

\begin{prop}
Assume that the base $S$ is a Noetherian scheme over $k$ 
and $Z \subset X$, $E$ are as before. 
For any point $s \in S$ let $F_s$ be the canonical extension 
of $E|_{X_s}$ to a locally free sheaf on $X_s$ and set 
$f(s) = \deg c_2(F_s)$.
Then $f: S \to \mathbb{Z}$ is a lower-semicontinuous function:
if $s_1$ is in the closure of $s_2$ then $f(s_1) \leq f(s_2)$. 
\end{prop}
\textit{Proof.} 
We can represent $E^\vee$ as a cokernel of the map 
$G_1|_U  \to G_0|_U$ where $G_i$ are 
direct sums of line bundles on $X$ (with ample dual line bundles).
Therefore on $U$ we have 
$$ 
0 \to E \to G_0^\vee \to G_1^\vee 
$$
Hence if $F$ is the direct image of $E$ to $X$ then 
$$
H^0(X, F) = Ker (H^0(X, G_0^\vee) \to H^0(X, G_1^\vee))
$$
Moreover, the same will hold for any base change $T \to S$: 
the cohomology of
the direct image extension of $E_T$ is always the kernel of
the morphism on global sections of $(G_0)_T^\vee$ and
$(G_1)_T^\vee$.  The same holds if we twist all sheaves by a 
relatively ample line bundle $L$ on $X$.  Since the rank of the morphism 
$\pi_* G_0^\vee (L) \to \pi_* G_1^\vee(L)$ is a lower 
semi-continuous function on $S$ (assuming both sheaves are 
locally free) we conclude that $h^0(F_s(L|_{X_s}))$ is an 
upper semi-continuous function on $S$. 
For an ample enough $L$ this is equal to the Euler characteristic 
of $F_s(L|_{X_s})$ which by Riemann-Roch is of the form 
$A - c_2(F_s)$ where $A$ is independent of $s$. Therefore 
$f(s) = c_2(F_s)$ is lower semi-continuous.
$\square$

\bigskip
\noindent
We see from the previous proposition that just deforming $E$
would give a badly behaved functor: a bundle with $c_2 = n$
may be deformed to a bundle with $c_2 = n+k$ and apriori 
there is no bound on the non-negative integer $k$.
Formally, this is manifested in the fact that $H^1(U, End(E))$
has an infinite dimensional piece coming from the local cohomology of
$End(E)$ at $x \in Z$. 
  However if
we know that $N$ is the product of $\mathfrak{c}_2^E$ and the
functor of an effective cycle of degree $k_0$ then it follows that
$k \leq k_0$. This keeps the tangent vectors of first order deformations
within a smaller subspace which is closely related to the standard
increasing filtration on the local cohomology $\mathcal{H}^2_x (End(E))$.

\subsection{Relation to (quasi)maps into a Quot scheme.}

Let us start with a situation when $E$ is defined everywhere on $X$,
choose a rank $r$ bundle $W$ on $X$ which is trivial on the fibers of 
$\pi$ (i.e. $W$ is a pullback of a bundle on $S$) and suppose that we have
an injective morphism of coherent sheaves $\rho: W \to E$ which is an isomorphism at the 
generic point of each fiber $X_s$,  $s \in S$. Under this assumption both
$\rho$ and its adjoint are embeddings of coherent sheaves on each fiber and
the cokernel $K := W^\vee/ E^\vee$ is flat a Cohen-Macaulay of pure relative
dimension 1 over $S$. So we can view $E^\vee$ as a kernel of the surjection 
$W^\vee \to K$ and using  Cramer's Rule we see that 
$K$ is annihilated by the section $\det (\rho)$ of 
$\det(E) \otimes \det(W)^\vee$, i.e. $K$ is supported on the 
zero scheme $Y$ of $\det(\rho)$. Hence $E^\vee$ is the kernel
of a composition 
$$
W^\vee \to  W^\vee |_Y \to K.
$$
Suppose that $L$ is sufficiently ample and set $H = \pi^* \pi_* L$, 
$P = \mathbb{P}(H)$. We view the
 incidence subscheme
$\mathcal{C} \subset P \times_S X$ as a relative curve $\mathcal{C} \to P$.
Consider the open
subset $\mathcal{U} \subset P$ of those curves that have no common
irreducible components with $Y$ (if $L$ is ample enough the closed 
complement to $\mathcal{U}$ will have high codimension).
Pull $K$ back to $P \times_S X$ and let $K|_{\mathcal{C}}$ be the restriction 
to $\mathcal{C}$. 
Then  $K|_{\mathcal{C}}$ is flat of finite length $l = \det(E) \cdot L$ 
over $\mathcal{U}$, since the restriction of $E^\vee$ to the fibers of $\mathcal{C}$
over $\mathcal{U}$ has degree $-l$. 
For $p \in \mathcal{U}$ and
 the corresponding curve $C_p \subset X$  we have a quotient map 
\begin{equation}
\label{quot}
W^\vee|_{Y \cap C_p}   \to K|_{C_p}.
\end{equation}
The source and target have length $r l$ and $l$, 
respectively, over the residue field of $p$.
Let $\mathcal{W}$ be the rank $rl$ vector bundle 
$
\pi_*( W^\vee|_{Y \cap \mathcal{C}})
$
on $\mathcal{U}$. The above construction defines a section 
$s: \mathcal{U} \to Gr(\mathcal{W}|_\mathcal{U}, l)$ of the 
Grassmannian bundle of rank $l$ quotients over $\mathcal{U}$. 

\bigskip
\noindent
Let $\mathcal{O}_{Gr}(1)$ be the determinant of the universal quotient 
bundle on the relative Grassmannian and $\mathcal{K}$ its pullback to 
$\mathcal{U}$ via $s$. Then $s$ can also be recovered from 
the induced surjection on $\mathcal{U}$
$$
\Lambda s: \Lambda^l \mathcal{W} \to \mathcal{K}
$$
which satisfies Pl\"ucker conditions. 
If we assume that $\Lambda^l \mathcal{W} \to \mathcal{K}$
is only generically surjective on fibers of $\mathcal{U} \to S$ (but still satisfies 
Pl\"ucker conditions), we get a concept of 
a \textit{quasimap} from $\mathcal{U}$ to the relative Grassmannian.

Usually, the concept of a quasimap is used when $P$ has relative dimension 
1 over $S$, i.e. it is a family of projective lines. In this setting $P$ would be a
projectivization not for the full bundle of sections of $L$, but its rank 2 
subbundle. 
Assume for simplicity that 
$S= Spec(k)$ and the 2-dimensional subspace $H \subset H^0(X, L)$ is chosen
so that $\mathcal{U} = P \simeq \mathbb{P}^1$ (i.e. we have a pencil of 
curves none of which shares a component with $Y$). Then, for a quasimap 
from $P$ to the Grassmannian $Gr(\mathcal{W}, r)$ the image of the morphism
$\Lambda s: \Lambda^l \mathcal{W} \to \mathcal{K}$ is a rank one locally 
free subsheaf of $\mathcal{K}$ with quotient supported at some points 
$p_1, \ldots, p_k \in \mathbb{P}^1$ and multiplicities $a_1, \ldots, a_k$.
 In other words, we can view
a quasimap $s$ as the usual map with pullback line bundle $\mathcal{K} (- \sum a_i p_i)$
which is ``twisted" by $\sum a_i p_i$, the defect of our quasimap.

\bigskip
\noindent
The point of this discussion is that, similarly to the case when a bundle $E$ defines
a map to the relative Grassmannian, a quadruple $(Z, E, N, D)$ with  an effective
$E$-localization defines a quasimap $\mathcal{U} \to Gr(\mathcal{W}|_{\mathcal{U}}, r)$.
The ``effectiveness" condition is responsible for the fact that $\Lambda s$ is a 
regular morphism of coherent sheaves, not rational. 
The main tool here is the following

\begin{lemma} If $E$ is defined everywhere on $X$ and $d = \deg c_2(E)$ then 
the line bundle $\mathcal{K}$ on $\mathcal{U} \subset P$ is isomorphic to 
$$
\det (\pi_{\mathcal{U}})_* (\mathcal{O}_{Y \cap \mathcal{C}})
\otimes \mathfrak{c}_2^E (L) \otimes \mathcal{O}_P(d)|_{\mathcal{U}}
$$
\end{lemma}
\textit{Proof} We will pull the bundle $\mathcal{K}$ back to 
the product $P \times_S G$ where $G = Gr(r-1, \mathcal{W})$ is the 
Grassmannian of $(r-1)$-dimensional subspaces in $\mathcal{W}$, and will show
a similar decomposition on that product, which will imply the decomposition
on $P$. The isomorphism to be proved  is based on the following diagram
of sheaves on $X \times G$ with exact columns ($Q$ being the dual to the universal 
rank $(r-1)$ sub-bundle on $G$ and $\mathcal{O}_G(1)$ the universal 
quotient bundle):
$$
\begin{CD}
0 @>>>\det( E^\vee) \otimes \mathcal{O}_G(-1)  @>>> \det(W^\vee) \otimes_k 
\mathcal{O}_G(-1)  @>>> K  \\
@. @VVV @VVV  @|\\
0 @>>> E^\vee  @>>> W^\vee @>>> K @ >>> 0 \\
@. @VVV @VVV @VVV\\
0 @>>> Q @= Q @>>> 0 \\
@. @VVV @VVV\\
@. \mathcal{A} @>>>  0 \\
\end{CD}
$$
Applying the Snake Lemma to the 
 two middle rows of the diagram,  we get a short exact sequence
$$
0 \to \mathcal{O}_Y \boxtimes \mathcal{O}_G(-1)  \to K 
\to \mathcal{A} \to 0.
$$
If we now pull this back to $X \times G \times P$, intersect with 
the preimage of the universal curve $C \subset X \times P$
and then push forward to $G \times P$, this will split the 
determinant for the direct image of $K$, into the tensor product
of $\Lambda^l \pi_* (\mathcal{O}_{Y \cap C}) \otimes
\mathcal{O}_G(-l)$ and
$\mathfrak{c}_2^{E^\vee} (L \boxtimes \mathcal{O}_P(1))
\otimes \mathcal{O}_G(l)$. Now the
assertion follows from a canonical isomorphism 
$$
\mathfrak{c}_2^{E} \simeq \mathfrak{c}_2^{E^\vee}
$$
(which is a well-known statement for cohomology classes, and for 
the multiplicative functors can be derived from Elkik's construction in 
a straightforward way).
$\square$

\bigskip
\noindent
This can be rephrased as follows (assume $S= Spec(k)$, $\mathcal{U}=P$
 to simplify things).
For $p$ in an open subset $\mathcal{U}_Z \subset 
P$ of curves avoiding $Z$ we have the same picture as before: 
a  length $l$ quotient $K_p$ of $W^\vee|_{C_p \cap Y}$ inducing a map 
$$
\mathcal{U}_Z \to Quot(W^\vee|_{C_p \cap Y}, l) \subset  
Gr(\mathcal{W}, l) \subset \mathbb{P}(\Lambda^l \mathcal{W})^\vee.
$$
So on $\mathcal{U}_Z$ we have a surjective morphism of vector bundles
$$
\Lambda^l \mathcal{W} \to \mathcal{K}.
$$
By the above, its target $\mathcal{K}$ extends to
a line bundle on  the whole $P$ as
$$
\det R\pi_*(\mathcal{O}_{C} \otimes^L  \mathcal{O}_Y) 
\otimes N(L) \otimes \mathcal{O}_P(d)
$$
and in fact $\Lambda^l \mathcal{W} \to \mathcal{K}$ extends to $P$ 
as well, although
the extension may not be surjective. This means that the map 
of $\mathcal{U}_Z$ to the relative Grassmannian 
$Gr(\mathcal{W}, l)$ is now extended to a \textit{quasimap}
on $P$, in the sense of \cite{FGK}. This provides 
a link with \cite{FGK} and \cite{BFG} since our Quot scheme 
naturally embeds into a locus of the affine Grassmanian 
of the curve $C_p$, where the elementary modifications of bundles 
are performed only at the points of $C_p \cap Y$. 

\section{The Uhlenbeck (or quasibundle) functor}

\subsection{Definition.}

Let $\pi: X \to S$ a smooth projective morphism of relative dimension 
2, as before. We define the groupoid category $QBun(r, d)$
of degree $d$ rank $r$
\textit{quasibundles} on $X$ over $S$, as follows. An object of this category
consists of  
\begin{enumerate}
\item A closed subset $Z \subset X$ which is finite over $S$ and
a rank $r$ vector bundle $E$ on the open complement $U = X \setminus Z$;

\item A line bundle $D$ on $X$ with a given isomorphism $D|_U 
\simeq \det(E)$;

\item A degree $d$ multiplicative functor $N: PIC(X/S) \to PIC(S)$ with an 
effective $E$-localization at $Z$.
\end{enumerate}
For two such objects, a morphism is given by an isomorphism 
$\psi: E_1 \to E_2$ over an open subset  $U^\circ \subset 
U_1 \cap U_2$, with closed
complement $Z^\circ$ finite over $S$, and an isomorphism of functors
$N_\psi: N_1 \simeq N_2$, which agrees with multiplicative 
structures and localizations at $Z^\circ$.  

Any base change $T \to S$ induces a quasibundle 
on $T$: for $N$ extension to base changes
 is a part of the definition of a multiplicative functor, 
and for $E, Z, D$ we can apply the natural base change definitions. 
This easily extends to the case when $S$ is an algebraic space 
(then, of course we use etale topology on $S$ instead of Zariski),
which is a bit more appropriate from the point of view of algebraic stacks.

\subsection{Completeness and points over $Spec(k)$.}

\begin{prop}
The functor $Bun(r, d) \subset QBun (r, d)$ of
rank $r$ bundles with $c_2 = d$, is open in $QBun(r, d)$.
\end{prop}
It follows from the definitions that $Bun(r, d)$ is formed by those
quasibundles for which $Z$ is empty and $N = \mathfrak{c}_2^E$. 

\bigskip
\noindent
\textit{Proof.} Let $(Z, E, N, D)$ be a quasibundle with base $S$
and suppose that after a base change $T \to U$ the pullback 
quasibundle $(Z_T, E_T, N_T, D_T)$ is isomorphic to a usual 
bundle. We need to show that $T$ factors through an 
open subscheme $S^\circ \subset S$ (possibly empty). 

By Section 4.4 for any point $s \in S$ the restriction of $E$ to 
the fiber $X_s$ extends canonically to a vector bundle $F_s$ with 
$c_2 = d -k$, $k \geq 0$. Moreover, the subset $S^\circ \subset S$
of all points for which $k =0$ is open in $S$. If the base change for
a quasibundle under $T \to S$ is isomorphic to an honest bundle, then
$E_T$ extends to a vector bundle with $c_2 = d$ along the fibers. 
Choosing a point $t \in T$ and restricting to its fiber, we see that 
the image of $t$ must be in $S^\circ$, hence the base change morphism
in fact factors through $T \to S^\circ$. 

In the opposite direction, let us assume that $S^\circ = S$ and prove
that the quasibundle is in fact an honest bundle. By Section 4.4 it 
suffices to show that the direct image $F$ of $E$ to $X$ is locally free. 
As in \textit{loc. cit.} after an ample twist of $E$ we have a sequence
$$
0 \to E \to G_0^\vee|_U \to G_1^\vee|_U
$$
where the two bundles on the left are direct sums of ample line bundles. 
Our assumption $S^\circ = S$ and the argument of Section 4.4 
imply that $\pi_* G_0^\vee\to \pi_* G_1^\vee$ has constant rank
hence $\pi_* F$ is locally free on $S$. This also holds for any ample
twist of $E$ and thus $F$ is flat over $S$. Then  the 
direct image extension
from $U$ to $X$. Since on each fiber $X_s$ such an extension $F_s$
its locally free, for $i \geq 1$ and $x \in X_s$ we have 
$\mathcal{E}xt^i_{X_s} (F|_{X_s}, k(x)) = 0$ where $k(x)$ is the skyscraper sheaf
at $x$. By flatness of $F$ and the change of rings spectral sequence 
 we get that $\mathcal{E}xt^i_X(F, k(x)) = 0$ for any $x \in X$.
Hence $F$ is locally free.
$\square$

\bigskip
\noindent
\textbf{Remark.} In fact, if $S_k \subset S$ is the locally closed subscheme
on which $c_2(F_s) = d - k$ and a base change morphism 
$T \to S$ factor through $S_k$ then one can show by a similar technique
that $E_T$ extends by direct image to a vector bundle $F$ with 
$c_2 = d-k$ along the fibers and the functor $N$ splits into the
functor $\mathfrak{c}_2^F$ and an functor induced by a relative
family of degree $k$ cycles on $T$. 

\bigskip
\noindent
The next result explains in which sense $QBun(r, d)$ is a compactification
of $Bun(r, d)$. Recall that for schemes the valuative criterion of 
properness has two parts: a certain morphism should admit an extension, 
and such extension should be unique. For functors of bundles 
uniqueness is too much to ask: a family of bundles over a punctured
curve has in general several non-equivalent extensions to the puncture.
On the other hand, one may still ask about existence of such an
extension. 

\begin{prop} If $\mathbb{O}$ is a discrete valuation ring  morphism $Spec(\mathbb{O}) \to S$ and $\mathbb{K}$ is its field of fractions
then each $\mathbb{K}$-valued point  of $QBun (r, d)$
lifts (non-uniquely!) to a $\mathbb{O}$-valued point of $QBun (r, d)$. 
\end{prop}
\textit{Proof.} A vector bundle over an open subset 
of $X_{\mathbb{K}} = X \times_S Spec(\mathbb{K})$
extends to a vector bundle $F$ on $X_{\mathbb{K}}$ itself by 
direct image. Hence $N$ will be of the form $\mathfrak{c}_2^F 
\otimes N_\xi$, for an effective cycle $\xi$ of certain degree $k$.
Since $\Gamma^k(X/S)$ is represented by a proper scheme 
over $S$, the cycle $\xi$ admits a unique extension over $\mathbb{O}$.
On the other hand, by Proposition 6 in \cite{La}  we can find a torsion
free sheaf $F_{\mathbb{O}}$, flat over $\mathbb{O}$, which 
restricts to $F_{\mathbb{K}}$ over $\mathbb{K}$. It will be 
locally free away from a closed subset of 
$X_{\mathbb{O}}$ which is finite over $Spec(\mathbb{O})$.
Also, the functor $\mathfrak{c}_2^{F_{\mathbb{O}}}$ is
well-defined since we can choose a locally free resolution
$0 \to F_1 \to F_0 \to F_{\mathbb{O}} \to 0$ and 
use the corresponding functors for $F_1, F_0$ (see next subsection
for details). Then the restriction of $F_{\mathbb{O}}$ 
and $N = \mathfrak{c}_2^{F_\mathbb{O}} \otimes N_{\xi_{\mathbb{O}}}$
extends our quasibundle to $Spec(\mathbb{O})$. $\square$

\begin{prop} If $\mathbb{K}$ is a field and $Spec(\mathbb{K}) 
\to S$ is a morphism of schemes then the set of isomorphism 
classes $|QBun(r, n) (Spec(\mathbb{K}))|$ has a set-theoretic
decomposition
$$
|QBun(r, d)(Spec(\mathbb{K}))| =
\coprod_{d' \geq 0} |Bun(r, d-d')(Spec(\mathbb{K}))|
 \times
\Gamma^{d'}(X)(Spec(\mathbb{K}))
$$
\end{prop}
\textit{Proof.}
A quasibundle over $Spec(\mathbb{K})$ always has the property 
that the vector bundle $E$ on $U_{\mathbb{K}} 
\subset X_{\mathbb{K}}$ has a locally free
direct image $F$. Applying Section 4.4 we see that a quasibundle in
this case is simply a pair $(F, \xi)$ where $\xi$ is an effective cycle
on $X_{\mathbb{K}}$ and $F$ is a vector bundle. 
Since $\deg c_2(F) + \deg \xi = d$, the assertion follows. 
$\square$

\subsection{Conjectured local covering by a scheme.} 

In this section we construct a scheme, or rather a system of schemes,
 mapping into the functor $QBun(r, d)$
and we conjecture that these give a faithfully flat covering 
for an increasing system of open substacks covering $QBun(r, d)$ (which would 
imply the $QBun(r, d)$ is an Artin Stack). 
We define an \textit{enhanced quasibundle} as the following data 

\begin{itemize}

\item A quasibundle $(Z, E, N, D)$ as defined before; 

\item A section $s$ of the line bundle $D$ extending $\det(E)$ from $U$ to $X$ with zero scheme 
$Y \subset X$ being of pure relative dimension 1 over $S$ and such that $Z \subset Y$; 

\item A line bundle $L$ together with a finite dimensional subspace
of section $H \subset H^0(X, L)$ which generate $L$. The projective 
space $P = \mathbb{P}(H)$ of lines in $H$ has two open subsets:
the subset $\mathcal{U}$ of sections of $L$ which are nonzero at 
the generic points of irreducible components of $Y$, and a smaller
open subset $\mathcal{U}_Z$ of those sections which are, in addition,  
nonzero at $Z$;

\item A vector bundle $W$ on $X$ pulled back from $S$ and a quasimap 
$\mathcal{U} \to Quot_{Y \cap C} (W^\vee|_{Y \cap C}, l)$ where $l= L \cdot Y$
 (the quasimap is understood in the 
sense of embedding the determinantal bundle on the $Quot$ scheme). 
We require  that the line bundle 
on $\mathcal{U}$ in the definition of a quasimap is 
$$
\Lambda^l \pi_* \mathcal{O}_{Y \cap C} \otimes N(L) \otimes \mathcal{O}_P(d) 
$$

\item Over the open subset $\mathcal{U}_Z$, we require 
that the quasimap is a morphism (i.e. a map in the usual sense) and that we
are given an 
 the isomorphism of 
$E^\vee$ with the kernel of the composition $W^\vee  \to W^\vee_Y \to Q$ where 
$Q$ is the quotient sheaf corresponding to the 
morphism from $\mathcal{U}_Z$ to the relative $Quot$ scheme. 

\end{itemize}

For a fixed $W$ and a relatively ample bundle $B$ on $X$ we conjecture that 
the subfunctor $QBun(r, d)_s$ of quasibundles such that $E \otimes B^{\otimes s}$
admits an enhancement for a choice of $s \geq 1$, form an increasing sequence of open 
subfunctors with the union equal to $QBun(r, d)$; and that the 
corresponding functors of enhanced quasibundles are represented by
schemes which provide a faithfully flat cover of $QBun(r, d)_s$ for $s \geq 1$.

\subsection{The Gieseker-to-Uhlenbeck morphism}

Observe that our construction of the functor $\mathfrak{c}_2$
may be generalized a bit. Namely, suppose that $E$ is a
coherent sheaf, flat over $S$, which restricts to a torsion
free sheaf on each fiber. Then at least locally on $S$
 we have a length 
2 projective resolution $0 \to E_1 \to E_0 \to E \to 0$
in particular it is a perfect complex, so that $\det E$ is
well-defined. Now we can use the same formula \eqref{defineC2}.

This defines a morphism of functors $Gies(r, d) 
\to QBun(r, d)$ where the source is the Gieseker 
compactification of flat families of torsion free sheaves
with rank $r$ and $c_2 = d$. Indeed, for any such sheaf
$E$ we can find $Z \subset X$ which is finite over 
$S$, such that $E$ is locally free on the open 
complement $U = X \setminus Z$, and the 
multiplicative 
functor $N = \mathfrak{c}_2^E$ constructed above
admits a natural $E$-localization at $Z$. 

We need to show that this localization is effective. 
The latter condition can be checked on the fibers of
$\pi$, i.e. we can assume that $S$ is the spectrum of 
some field extension of $k$. But then $E^{\vee \vee}$ 
is locally free with $A = E^{\vee \vee}/E$ an Artin sheaf.
Since 
$$
\mathfrak{c}_2^E \simeq\mathfrak{c}_2^{E^{\vee \vee}}
\otimes N_A,
$$
the assertion follows.

\section{Simply-connected semi-simple groups.}

\subsection{Principal quasibundles on surfaces.}

Let $G$ is split semi-simple and connected group over $k$ (automatically
split since $k$ is closed), with a
maximal torus $T$,  Weyl group $W$, character lattice
$\mathbb{X}$ and the dual lattice $\mathbb{Y}$.
Denote also by $Q(G)$ the lattice of $W$-invariant symmetric bilinear
forms on $\mathbb{Y}$. Recall that in our case such forms are automatically 
even and that $Q(G)$ is a free abelian group of rank $t$ 
equal to the number of almost simple factors of $G$. 
 We denote by 
$\pi(G)$ the dual free abelian group
(it is isomorphic to the third homotopy group of
$G$ for $k = \mathbb{C}$). Both groups do not depend on the 
choice of $T$ and a homomorphism 
$\rho: G_1 \to G_2$ of split semisimple simply connected groups
induces a homomorphism $\rho_*: \pi(G_1) 
\to \pi(G_2)$.  The meaning of $Q(G)$ is clarified by the following result, 
cf. Theorem 4.7 in \cite{BD}:
\begin{prop}
Under the assumptions above $Q(Y) \simeq  H^1_{Zar}(G, \mathcal{K}_2)$. 
Any $\mathcal{K}_2$-torsor on $G$ has a unique multiplicative structure and 
no nontrivial automorphisms. 
\end{prop}
Observe also that the group $Q(G)$ contains a semigroup
of non-negative quadratic forms, isomorphic to $(\mathbb{Z}^{\geq 0})^t$,
and consequently $\pi(G)$ also has the dual semigroup $\pi_+(G)$ of 
non-negative elements. By the above proposition, on $G$ we have a universal 
torsor over $\mathcal{K}_2 \otimes_{\mathbb{Z}} \pi(G)$ which we denote
by $\mathfrak{C}_2$.

\begin{prop}
A choice of a principal $G$-bundle $P$ on $X$ in the Zariski topology induces
a multiplicative functor 
$$
\mathfrak{c}_2^P: PIC(X/S) \to PIC(S) \otimes_{\mathbb{Z}} \pi(G)
$$
where the target category is understood as torsors over 
$\mathcal{O}^* \otimes_{\mathbb{Z}} \pi(G) \simeq 
(\mathcal{O}^*)^{\times t}$. It also induces 
a locally constant degree function $d_P: S \to \pi(G)$ (the generalized second Chern class of $P$).
For a group homomorphism $\rho: G_1 \to G_2$ 
of semi-simple simply-connected groups one has 
a natural isomorphism of functors
$$
 \rho_* \circ \mathfrak{c}_2^P \simeq \mathfrak{c}_2^{\rho(P)}
$$
where $\rho(P)$ is the induced principal bundle over $G_2$. 
Also, $d_{\rho(P)} = \rho_* \circ d_P$. 
\end{prop}
\textbf{Remark.} 
At this moment we cannot yet prove the key property 
that  for a line bundle $L$ on $S$ one has 
$$
\mathfrak{c}_2^P(\pi^* L) \simeq L^{\otimes d_P}.
$$
which is to say,  that $\mathfrak{c}_2^P$ is a degree $d_P$
multiplicative functor.

\bigskip
\noindent
\textit{Proof.} It suffices to choose an element $\alpha$
in the dual group $Q(G)$ and construct a function 
$S \to \mathbb{Z}$ and a functor $PIC(X/S) \to PIC(S)$ which 
is additive in $\alpha$. 

By the work of Brylinski and Deligne,  cf. \cite{BD}, 
such a choice of $\alpha$
gives a central extension $\widetilde{G}$
of $G$ by $\mathcal{K}_2$, and hence 
an additional choice of a $G$-torsor $P$ on $X$ 
induces a 
a $\mathcal{K}_2$-gerbe 
$\mathfrak{C}_2(P, \alpha)$ on $X$ with a class in 
$H^2(X, \mathcal{K}_2)$ which we can call the second Chern 
class of $P$ corresponding to $\alpha$. This gerbe can either be  understood as the gerbe of lifts
of $P$ to a torsor over $\widetilde{G}$; or we can consider the canonical ``cocycle map"
$P \times_X P \to G$ and then pull back 
the $\mathcal{K}_2$-torsor of local 
splittings $G \to \widetilde{G}$. The result 
will be a ``bundle gerbe", cf. \cite{St}
with the product structure on 
$P \times_X P \times_X P$ induced by the 
product on $\widetilde{G}$. 

Now the function $\alpha \circ d_P: S \to \mathbb{Z}$ is 
defined easily. For any closed point $s \in S$ the
class of $\mathfrak{C}_2(P, \alpha)$  projects to $H^2(X_s,
\mathcal{K}_2)$ which is the Chow group of zero
cycles cycles up to rational equivalence. Then the function in question
simply computes the degree fiber by fiber.

Ideally, the construction of the functor would proceed as follows.
An object $L$ in $PIC(X/S)$, i.e. a torsor over $\mathcal{K}_1$
can be multiplied by $\mathfrak{C}_2(P, \alpha)$ via external 
product, to give a 2-gerbe over $\mathcal{K}_2 \otimes_Z \mathcal{K}_1$
and then the K-theory product induces a 2-gerbe over $\mathcal{K}_3$, 
and even over Milnor $K$-theory
(the 2-gerbe can be understood, e.g. as a bundle 2-gerbe in the sense of 
\cite{St}. Such an object has a class in $H^3(X, \mathcal{K}_3)$. 
Since $\pi: X \to S$ has two dimensional fibers and we are using 
Zariski topology, this class projects to 
$H^1(S, R^2 \pi_* \mathcal{K}_3)$. Using a slightly more refined
geometric language, locally on $S$ the above 2-gerbe over $\mathcal{K}_3$
can be trivialized (or neutralized) and all such trivializations form 
a torsor over $R^2\pi_* \mathcal{K}_3$. If we had a 
``residue-trace" map $R^2\pi_* \mathcal{K}_3 \to \mathcal{K}_1$
similar to integration of differential forms along the fibers of $\pi$, 
this would induce a $\mathcal{K}_1$-torsor on $S$, i.e. a line bundle. 

Since the ``residue-trace" map has not been constructed at the moment 
of writing, we choose a roundabout strategy. First, as in 
\cite{El} and \cite{MG} it suffices to define the functor on 
bundles $L$ which are ample enough, and extend to arbitrary bundles 
by multiplicativity. Hence we assume that $L$ is relatively very ample
with locally free $\pi_* L$ and vanishing higher direct images. 

Next, we replace $S$ by the projective bundle $P$ of lines in 
$\pi_*(L)$ and then consider the tautological curve $Y \subset X_P$. 
We will construct a certain bundle on $P$ which necessarily has 
the form $\xi^{\otimes k} \otimes \mu^* N$ where 
$\xi$ is the tautological relative line bundle for the morphism
$\mu: P \to S$ and $N$ is the a certain line bundle on $S$ which we
declare to be $\alpha \circ \mathfrak{c}_2^P(L)$. Exactly the same 
strategy has been used in \cite{El} to construct ``intersection bundles".

Consider the restriction of the gerbe $\mathfrak{C}_2(P, \alpha)$
to $Y$ and call this restriction $\mathfrak{C}_2$ for short. 
Since the fibers of the projection are one dimensional, the 
gerbe  $\mathfrak{C}_2$ can be trivialized locally on $P$ and such 
trivializations induce a torsor over $R^1 \eta_* \mathcal{K}_2$.
On the cohomological level, the term $H^0(S, R^2 \eta_* \mathcal{K}_2)$
in the Leray-Serre spectral sequence of $\eta$, vanishes. Hence the 
class of  $\mathfrak{C}_2$ in $H^2(Y, \mathcal{K}_2)$ projects to 
$H^1(S, R^1 \pi_* \mathcal{K}_2)$. 

Next we use the $K$-theory localization sequence
to construct a morphism $R^1 \eta_* \mathcal{K}_2 
\to \mathcal{K}_1$. Indeed, on $Y$ we have a long exact sequence of 
sheaves
$$
\ldots \to \mathcal{H}_i \to \mathcal{Q}_i  \to \mathcal{K}_i 
\to \mathcal{H}_{i-1} \to \ldots
$$
where $\mathcal{Q}_i$ are Zariski sheafified $K$-theory groups
for the exact category of projective modules over the localization 
$T^{-1} \mathcal{O}_Y$
of $\mathcal{O}_Y$ at the multiplicative set $T$
of all relative non-zero divisors
(i.e. regular functions on $Y$ which are nontrivial on any irreducible
component of a fiber over $P$). Similarly, $\mathcal{H}_i$
is the Zariski sheafified $i$-th $K$-theory for the exact category of sheaves
on $Y$ which are finite and flat over $P$ and admit a length 1 projective
resolution as coherent sheaves on $Y$. See \cite{Gr4} for details on localization of projective modules in $K$-theory. 

We are interested in the piece of the short exact sequence given by 
$$
\ldots \to \mathcal{K}_2  \to \mathcal{Q}_2 \to \mathcal{H}_1 
\to \ldots
$$
we don't know whether the dots can be replaced by zeros 
(which would follow from a relative version of Gersten conjecture) but
we can replace $\mathcal{K}_2$ by its image $\mathcal{A}_2$ in 
$\mathcal{Q}_2$ and $\mathcal{H}_1$ by the image $\mathcal{B}_1$
of $\mathcal{Q}_2$ and get a short exact sequence $0 \to \mathcal{A}_2
\to \mathcal{Q}_2 \to \mathcal{B}_2 \to 0$.

Next we note that the sheaf $\mathcal{Q}_2$ satisfies the relative
version of flasqueness: any section defined on a neighborhood $V$
of a point $y \in Y$ extends to the preimage of a small neighborhood
of $\eta(y) \in P$. Indeed, we are dealing with sheafified $K_2$
for a localization $T^{-1} \mathcal{O}_Y$ 
of $\mathcal{O}_Y$ and $K_2$ of a local ring may
be described via the Steinberg presentation, so the flasque property 
of $\mathcal{Q}_2$ follows from that of $T^{-1} \mathcal{O}_Y$. 

Therefore $R^1\eta_* \mathcal{Q}_2 = 0$ on $P$ and we have
$R^1 \eta_* \mathcal{K}_2 = Coker(\eta_* \mathcal{Q}_2 
\to \eta_* \mathcal{B}_1)$. Hence to construct a morphism 
$R^1 \eta_* \mathcal{K}_2 \to \mathcal{K}_1$ it suffices to 
construct a morphism $\eta_* \mathcal{H}_1 \to \mathcal{K}_1$
which vanishes on the image of $\eta_* \mathcal{Q}_2$ in 
$\eta_* \mathcal{B}_1 \subset \mathcal{H}_1$. 

But $\mathcal{K}_1$ was obtained by sheafifying $K_1$ of the
exact category of sheaves on $Y$ which are finite and flat over $P$ and 
have projective dimension 1, hence the usual direct image of a sheaf
does induce a sheaf morphism $\eta_* \mathcal{H}_1 \to \mathcal{K}_1$.
The fact that it vanishes on the image of sheafified $K_2$ for
$T^{-1} \mathcal{O}_Y$ follows from Section 7 of \cite{Gr2} which 
computes the corresponding connecting homomorphism of the localization 
sequence. 

The functorial behavior of $\mathfrak{c}_2^P$ and $d_P$ with 
respect to $\rho: G_1 \to G_2$ follows from the functorial behavior of 
the central extensions by $K_2$, as established in \cite{BD}. 
$\square$.

\bigskip
\noindent
\textbf{Remarks.}
(1) We remark here on the assumption of Zariski local triviality of $P$.
When $S = Spec(k)$, by the work of de Jong, He and Stuart, cf. \cite{dJHS}, any
etale  $G$ torsor $P$ has a rational section, and then by \cite{Ni} it is actually Zariski 
locally trivial. It is natural to conjecture that for arbitrary $S$ there is an etale 
covering  $T \to S$ such that the pullback $P_T = P\times_S T$  is
Zariski locally trivial. If this conjecture holds $\mathfrak{c}_2^P$ would be 
constructed by etale descent (even in the case of algebraic spaces).

Although the sheaves $\mathcal{K}_j^{p/(p+1)}$ are not flasque in 
general, by the construction of K-theory groups in terms of 
acyclic binary multicomplexes, cf. \cite{Gr3}, it is fairly easy to see that
these sheaves are indeed flasque over $S$, which gives the 
last feature that we need to construct the functor
$\mathfrak{c}_2^P$ and its degree function.

\bigskip
\noindent
(2) The same functor $\mathfrak{c}_2^P: PIC(X/S) 
\to PIC(S) \otimes_{\mathbb Z} \pi(G)$ may be interpreted as follows. 
Assume that $L$ is a line bundle on $X$ and its section $s$ 
 vanishes on a relative divisor $C \subset X$ which is  a 
\textit{smooth} family of curves over $S$. Then $P|_C$ defines a 
morphism of $S$ to the relative moduli stack of $G$-bundles on 
$C$, and since $\pi(G)$ is dual to the Picard group of this moduli 
stack, cf. \cite{BH}, there is a pullback $\mathcal{O}^* \otimes_{\mathbb{Z}} 
\pi(G) $ torsor on $S$.  In fact, we could attempt to \textit{define} our functor 
by choosing an appropriate open covering of $S$ and a section $s_i$ over each 
open set of the covering. But then one has to show independence of the choice 
of $s_i$. Of course, another choice $s'_i$ may be connected to $s_j$ by 
a copy of $\mathbb{P}^1$ in the space of sections, but some of the 
sections in the family connecting $s_i$ and $s_i'$ will not give smooth curves over
$S$. The best we can hope for is irreducible curves with nodal singularities. 
Since the corresponding torsors over the moduli of principal $G$ bundles 
on nodal curves have not been constructed yet, we had to resort to the  approach as above.

\bigskip
\noindent
Now the definiton of a $G$-quasibundle on $X$ follows a pattern
similar to the case of vector bundles. More exactly, for a given
locally constant degree function $d: S \to \pi(G)$
we consider 
a bundle $P$ on the open complement $U \subset X$ to a closed
subset $Z$ which is finite over $S$, and a degree $d$ 
multiplicative functor 
$$
N: PIC(X/S) \to PIC(S) \otimes_{\mathbb{Z}} \pi(G).
$$
As in the case of vector bundles, we want to consider a $P$-localization of 
$N$ at $Z$, i.e. for any section $s: \mathcal{O}_X \to L$
of a line bundle on $X$ with zero scheme
$C$  disjoint from $Z$, we want to fix an isomorphism 
$$
a_s: \mathfrak{c}_2^P (L) \simeq N(L)
$$
which behaves multiplicatively with respect to tensor products
of line bundles and their sections. The LHS is well defined
since $P$ exists in on open neighborhood of $C$. 

Next, observe that any representation $\rho: G \to SL(r)$
induces a vector bundle $E=\rho(P)$, and a multiplicative functor 
$\rho(N): PIC(X/S) \to PIC(S)$ of 
degree $\rho(d)$, which comes equipped with
a natural $E$-localization. 

\medskip
\noindent
\textbf{Definition.} We will say that a $P$-localization of $N$ is
\textit{effective} if the induced $\rho(P)$-localization of $\rho(N)$
is effective for any choice of a representation $\rho: G \to SL(r)$. 

\medskip
\noindent
\textbf{Remark.} An alternative definition of an effective localization 
may be phrased as follows. Suppose that, after increasing $Z$, we can 
find a $B$-bundle $P_0$ inducing $P$ on $X \setminus Z$ (where 
$B$ is a Borel subgroup of $G$). Suppose further that the $T$-bundle 
bundle $P_1$ induced from $P_0$ via $B \to B/U \simeq T$ does 
extend to a bundle on $X$ which we also denote by $P_1$. 

Write the ``universal symmetric $W$-invariant form" on $\mathbb{Y}$ with
values in $\pi(G)$ as an element of $\pi(G) \otimes_{\mathbb{Z}}
\mathbb{X} \otimes_{\mathbb{Z}} \mathbb{X}$ in a certain basis $x_1, x_2, \ldots$
of $\mathbb{X}$:
$$
\sum a_{ij} \otimes x_i \otimes x_j, \qquad a_{ij} \in \pi(G).
$$
Since each $x_i$ is a character on $T$ we can use it to produce a line bundle 
$L_i$ from $P_1$, which is defined everywhere on $X$. With this notation, set
$$
N_{P_1} = \prod IB_{L_i, L_j}^{\otimes a_{ij}}: PIC(X/S) \to PIC(S) \otimes_{\mathbb{Z}}
\pi(G).
$$
One can show that, for our original $P$-localized $N$, the functor $N \otimes N_{P_1}^{-1}$ is
$\mathcal{O}$-localized at $Z$, i.e. corresponds to a relative zero cycle over $S$ but
not necessarily effective. So we require that for every choice of a  $B$-structure on $P$ 
(perhaps after a base change on $S$ and shrinking of an open set $U$ on which $P$ is 
defined), the quotient of the two multiplicative functors \textit{is} given by an effective
relative zero cycle. However, this definition is useful only if: (a) one proves a higher dimensional
version of Drinfeld-Simpson lemma (which would state that after an etale base change a $G$-torsor 
admits a $B$-structure away from codimension two); (b) one has a useful technique for comparing
different $B$-structures on the same $G$.

\medskip
\noindent
\textbf{Definition.} For a semi-simple simply connected $G$ and a choice of $d$, a degree $d$ 
\textit{quasibundle over $G$} is defined as a triple $(Z, P, N)$ consisting of a closed subset
$Z \subset X$ which is finite over $S$, a principal $G$-bundle $P$ on $U = X \setminus Z$, 
 a degree $d$ multiplicative functor
$N: PIC(X/S) \to PIC(S) \otimes_{\mathbb{Z}} \pi_3(G)$; plus 
 an effective $P$-localization of $N$ at $Z$. We denote by 
$QBun(G, d)$ the functor of such quasibundles. 

\bigskip
\noindent
The introduced functor has the following properties.

\begin{prop}
The functor $QBun(G, d)$ has the following properties:
\begin{enumerate}
\item A group homomorphism $\rho: G_1 \to G_2$ induces 
a morphism of functors $QBun(G_1, d) \to
QBun(G_2, \rho_*(d))$;

\item $QBun(G_1 \times G_2,  d_1 \times d_2)
\simeq QBun(G_1, d_1) \times QBun(G_2, d_2)$;

\item Assume that for a $G$-bundle $P$ the Chern class functor
$\mathfrak{c}_2^P$ has degree $d_P$. 
For an algebraically closed field $\mathbb{K}$ with a morphism   $Spec(\mathbb{K}) \to S$, the set of 
isomorphism classes of $\mathbb{K}$-valued points
 has a set-theoretic decomposition 
$$
|QBun(G, d) | = \coprod_{l \in \pi_+(G)} 
|Bun(G, d-l)| \times \Gamma^l(X)
$$
\end{enumerate}
where for $l = (l_1, \ldots, l_t) \in 
(\mathbb{Z}^{\geq 0})^t \simeq\pi_+(G)$ we denote by 
$$
\Gamma^l(X) \simeq \Gamma^{l_1} (X) \times_S \ldots
\times_S \Gamma^{l_t}(X/S)
$$
the space of $\pi_+(G)$-effective $0$-cycles on $X$. 
\end{prop} 
\textit{Proof.}
Indeed $\rho$ gives an induction functor on torsors, composition on 
degree functions and multiplicative functors, which shows (1). Part (2)
is also straigtforward as well since all our definitions are compatible with 
direct products (e.g. a torsor over $G_1 \times G_2$ is a pair $(P_1, P_2)$
consisiting of a torsor over $G_1$ and a torsor over $P_2$; a function
$d: S \to \pi(G_1 \times G_2)$ is a pair of functions $(d_1, d_2)$
with values in $\pi(G_1)$ and $\pi(G_2)$, respectively. 
For part (3), observe that a $G$-torsor $P$ on 
$U \subset X_{\mathbb{K}}$
always extends to $X_{\mathbb{K}}$ for any semisimple (or even 
reductive $G$). This follows from the Tannakian description of $G$-bundles as tensor functors from $G$-representations to vector bundles
on $X_{\mathbb{K}}$ 
since every induced vector bundle of $P$ extends to the whole surface 
simply by direct image as coherent sheaf. Hence, for a quasi-bundle 
$(Z, P, N)$ on $X_{\mathbb{K}}$ we can compare the multiplicative
functor $N$ with the functor $\mathfrak{c}_2^P$. By 
 our assumption on the degree of 
$\mathfrak{c}_2^P$, we get a multiplicative functor
$N \otimes (\mathfrak{c}_2^P)^{-1}$ of degree $d - d_P$. 
By effectiveness of the localization and the main result of \cite{Ba2}
this must be given by an effective cycle with coefficients in 
$\pi(G)$ (in particular, $d - d_P$ automatically takes values in $\pi_+(G)$).
$\square$

\bigskip
\noindent
We also formulate the following 

\bigskip
\noindent
\textbf{Conjecture}
\begin{enumerate}

\item There
 exists an injective group homomorphism $\rho: G \to \prod_i SL(r_i)$ such that the induced
morphism of functors 
$$
QBun(G, d) \to \prod_i QBun(r_i, \rho_i(d))
$$ 
is a closed embedding;

\item The functor $Bun(G, d)$ of principal 
$G$-bundles with characteristic class $d$, is an open subfunctor of $QBun(G, d)$. The functor $QBun(G, d)$
is complete, i.e. for a DVR $\mathbb{O}$ with the
fraction field $\mathbb{K}$, a family of quasibundles
over $Spec(\mathbb{K})$ extends (non-uniquely) 
to a family of quasibundles over $Spec(\mathbb{O})$.

\end{enumerate}

\medskip
\noindent
\textbf{Remarks.}
(1) Following a referee's suggestion, we mention one delicate issue related to the above conjecture. For some simple 
$G$ even on numerical level it will not be true that $c_2(P)$ can 
be obtained as a linear combination of $c_2(E_i)$ for various 
vector bundles induced from $P$ via a collection of representations
$\rho_i$. In general we can recover only a multiple $m c_2(P)$ where
$m$ is the greatest common divisor of Dynkin indices of all 
representations (in general, greater than 1). Consequently, 
if $N$ is the multiplicative functor of a quasi-bundle and $N_i$
are induced multiplicative functors for vector bundles $E_i$, we
can only recover $N(L)^{\otimes m}$ from $N_i(L)$ alone. 
After that, on needs to extract an appropriate $m$-th root of that bundle
which may not exist for a general bundle, and it may not be unique if
it exists. We note, however, that at least non-uniqueness is not a problem 
since for relatively ample $L$ we can find sections with zero set 
avoiding $Z$, and get $N(L)$ without any ambiguity by $P$-localization. 
In somewhat different language, we are using the fact that the map
on zero cycles sending $\xi$ to $m\cdot \xi$ is a closed embedding 
$\Gamma^l(X) \to \Gamma^{m\cdot l}(X)$.

(2) Second part of the above conjecture follows from the first fairly 
easily since the corresponding facts are known for usual 
(i.e. non-principal) quasi-bundles

(3) In the case of general reductive $G$ the situaltion is more delicate.
See \cite{BD} for the description of the category of multiplicative
$\mathcal{K}_2$-torsors on $G$ and \cite{BH} for the related description 
of the Picard group of the moduli stack of $G$-bundles on curves. 
We only note here that 
in defining quasibundles over this general $G$ one needs to 
consider a $G$ bundle $P$ on the open subset $U$ as before, but 
also fix an extension to $X$ of the induced principal bundle
over $G/[G,G]$ (e.g. the determinant line bundle for $G = GL(r)$).

\subsection{Remarks on quasibundles in higher dimensions.}

Let $\pi: X \to S$ be a smooth projective morphism with
connected fibers of fixed dimension $e \geq 2$ (with some additional
effort one could assume that $X$ is only flat and projective over $S$, and 
satisfies Serre's conditon $S_2$).  
A quasibundle on $X$
is still defined as a bundle $P$ over an open subset $U$ with closed
complement $Z$ of codimension $\geq 2$ over $S$, plus an
appropriate multiplicative functor $N$ with an effective $P$-localization
at $Z$. The multiplicative functor $N$ is now of the form
$$
N: PIC(X/S)^{\times (e-1)} \to PIC(S) \otimes_{\mathbb{Z}} \pi(G)
$$
Its degree is defined as a homomorphism
$$
d: Pic(X/S)^{\times (e-2)} \to \underline{\pi(G)}
$$
where $Pic(X/S)$ stands for the set of isomorphism classes in $PIC(X/S)$
and $\underline{\pi(G)}$ is understood as the sheaf of locally 
constant functions with values in $\pi(G)$. More explicitly, 
we require that 
$$
N(\pi^*L, L_1, \ldots, L_{e-2}) \simeq
L^{\otimes d(L_1, \ldots, L_{e-2})}, 
$$
and that $N$ is symmetric in its arguments in the sense explained in \cite{Du}. 

\bigskip
\noindent
When $P$ is defined everywhere on $X$, we could attempt to define such a functor such a functor $\mathfrak{c}_2^P$ 
in two ways, but each has technical difficulties at the moment. First we can 
use the product in K theory and multiply $\mathcal{K}_1$-torsors $L_0, \ldots L_{e-2}$ 
by the $\mathcal{K}_2$-gerbe $\mathfrak{C}_2(P)$ to get an $e$-gerbe over 
$\mathcal{K}_{e+1}$, i.e. a geometric object having a class in $H^{e+1}(X, \mathcal{K}_{e+1})$.
On each closed fiber over $S$ such an $e$-gerbe can be trivialized, and 
all local trivializations induce an etale torsor over $S$ with the structure group
$R^{e} \pi_* \mathcal{K}_{e+1}$. Next, one needs a general relative Gersten resolution 
$$
0 \to \mathcal{K}_{e+1} \to \mathcal{K}^{0/1}_{e+1} \to 
\mathcal{K}^{1/2}_{e} \to \ldots \to \mathcal{K}^{e/(e+1)}_1 \to 0
$$
with the property that the sheaves $\mathcal{K}_j^{(e+1-j)/(e+2-j)}$ have vanishing
higher direct images $R^q\pi_*$, $q > 0$. In that case 
$R^{e} \pi_* \mathcal{K}_{e+1} \simeq \pi_* \mathcal{K}^{e/(e+1)}_1$ has a natural 
norm map to $\mathcal{K}_1$ (exactly as for $e=2$) and one gets a $\mathcal{K}_1$-torsor
on $S$, which is the image of our functor. Here the technical difficulty is obviously
constructing the relative Gersten resultion and proving the vanishing of higher direct images. 

Alternatively, one can try to choose sections of $L_0, \ldots, L_{e-2}$ (locally over $S$) 
such that its common zero scheme $C \subset X$ is a family of smooth curves over $S$.
Then $S$ will map to the stack of $G$-bundles on the relative curve, and our 
functor is the pullback of the standard bundle on that stack. Here, as for the $e=2$ case
the difficulty is to show independence of the choice of sections, and this needs  
existence of the required bundle for a family of curves with possible nodal singularities. 

\bigskip
\noindent
A general $G$-quasibundle will be, as before, a triple $(Z, P, N)$ where 
$Z \subset X$ is a closed subset of pure relative dimension $(e-2)$ over $S$, 
$P$ a $G$-bundle over $U = X\setminus Z$, and $N$ a multiplicative functor
with $(e-1)$ arguments, as above, \textit{plus} an effective $P$-localization of
$N$ at $Z$, an identification of $N(L)$ and $\mathfrak{c}_2^P(L)$ for any 
section of $L$ which is nonzero along $Z$. Effectiveness for principal bundles 
can be defined by looking at all induced vector bundles, and for a vector bundle
$E$ one requires that each short exact sequence $0 \to \mathcal{O}_X^{(r-1)}
\to E \to \det (E) \to 0$ (which is chosen after a possible ample twist
of $E$, base change on $S$ and enlargement of $Z$), identifies the corresponding 
twist of $N$ with the multiplicative functor of a codimension 2 \textit{effective} cycle
supported at $Z$. This depends on the fact that relative effective cycles can be 
described via similar multiplicative functors.

The product of zero cycle spaces $\Gamma^d(X/S)$
is now replaced by the Chow \textit{scheme} $Chow^d(X/S)$ 
of multiplicative degree $d$ functors with an 
appropriate $\mathcal{O}$-localization (which is a correct definition for the
family of effective cycles in arbitrary characteristic). 
We expect that the functor $QBun(G, d)$ has a
completeness property (i.e. a family over a punctured spectrum of a 
DVR extends - non-uniquely - to the puncture), and a 
decomposition for the set of points over $Spec(k)$:
$$
QBun(G, d) = \coprod_{k \in \pi(G)} Bun'(G, d-k)
\times Chow^k(X/S)
$$
where $Bun'(G, \cdot)$ is the functor of $G$-bundles defined on 
an open subset of $X$ with closed complement of codimension 
$\geq 3$. The latter is an Artin stack by \cite{Ba2}.

\end{document}